\documentclass[11pt,twoside]{article}
\usepackage{amsmath,amssymb} 
\usepackage{amsfonts}
\usepackage{xcolor}
\usepackage[top=2.2cm, bottom=0.8in, left=2.25cm,
right=2.25cm]{geometry}
\usepackage{mathtools}
\usepackage{tikz-cd}
\usepackage{stmaryrd}

\usetikzlibrary{positioning}
\usetikzlibrary{arrows}

\newcommand{\R}{\mathbb{R}}

\newcommand{\X}{\mathfrak{X}}

\DeclareMathOperator{\id}{id}

\usepackage{amsthm}     
\newtheorem{theorem}{Theorem}
\newtheorem{lemma}[theorem]{Lemma}
\newtheorem{proposition}[theorem]{Proposition}

\theoremstyle{definition}

\newtheorem{example}[theorem]{Example}
\newtheorem{definition}[theorem]{Definition}

\newcommand{\la}{\lambda}

\newcommand{\ph}{\varphi}

\newcommand{\Ga}{\Gamma}

\DeclareMathOperator*{\equals}{=}
\DeclareMathOperator*{\addition}{+}
\DeclareMathOperator*{\subtraction}{---}
\DeclareMathOperator*{\mult}{\cdot}

\DeclareMathOperator*{\comma}{,}

\newcommand{\act}{\mathbin{\hbox{$<\kern-.4em\mapstochar\kern.4em$}}}
\newcommand{\ract}{\mathbin{\hbox{$\mapstochar\kern-.3em>$}}}

\newcommand{\dual}{(*)}
\newcommand{\ced}[1]{\triangleright #1}  
\newcommand{\bolt}{\text{\lightning}} 
\newcommand{\leftbracket}{[}
\newcommand{\rightbracket}{]}
\newcommand{\hlift}{H}
\newcommand{\wtilde}{\widetilde}
\DeclareMathOperator{\warp}{w}           
\DeclareMathOperator{\uwarp}{u}  
\newcommand{\lpair}{\llbracket}
\newcommand{\rpair}{\rrbracket}         

\usepackage{calrsfs}
\DeclareMathAlphabet{\pazocal}{OMS}{zplm}{m}{n}

\definecolor{kblue}{rgb}{0.2,0.2,0.6}
\definecolor{mpink}{rgb}{0.9,0.67,0.7}
\definecolor{mpink2}{rgb}{0.8,0.57,0.6}

\usepackage{xspace}

\newcommand{\tvb}{triple vector bundle\xspace}
\newcommand{\tvbs}{triple vector bundles\xspace}
\newcommand{\dvb}{double vector bundle\xspace}
\newcommand{\dvbs}{double vector bundles\xspace}
\newcommand{\vb}{vector bundle\xspace}
\newcommand{\vbs}{vector bundles\xspace}

\newcommand{\Tvbs}{Triple vector bundles\xspace}

\newcommand{\baf}{\mathsf{BF}} 
\newcommand{\lr}{\mathsf{LR}}  
\newcommand{\ud}{\mathsf{UD}}  

\usepackage{wasysym}
\newcommand{\duer}{\mathbin{\raisebox{3pt}{\varhexstar}\kern-3.70pt{\rule{0.15pt}{4pt}}}\,}

\usepackage{fancyhdr}

\renewcommand{\Bar}[1]{\overline{#1}}

\newcommand{\one}{\mathsf{ZYX}}
\newcommand{\two}{\mathsf{YZX}}
\newcommand{\thr}{\mathsf{XZY}}
\newcommand{\fou}{\mathsf{ZXY}} 
\newcommand{\fiv}{\mathsf{YXZ}}
\newcommand{\six}{\mathsf{XYZ}}

\usepackage{subfig}
\usepackage[labelsep=period]{caption}[2006/03/21]

\newcommand{\co}{\colon\thinspace} 

\newcommand{\plus}[1]{\addition\limits_{#1}}

\newcommand{\gog}{\mathfrak{g}}

\newcommand{\llangle}{\langle\!\langle}
\newcommand{\rrangle}{\rangle\!\rangle}

\usepackage{url}

\begin{document}

\title{\textbf{Warps and grids for double and triple vector bundles}\footnote{Mathematics 
Subject Classification (MSC2000): 53C05 (primary), 18D05, 18D35, 55R65 (secondary).}}

\author{Magdalini K. Flari and Kirill Mackenzie\\[2mm]
School of Mathematics and Statistics\\
University of Sheffield\\
Sheffield, S3 7RH, United Kingdom
}

\date{}

\maketitle

\vspace*{-12mm}
\begin{center}
\url{mkflari1@sheffield.ac.uk} \\
\url{K.Mackenzie@sheffield.ac.uk}
\end{center}

\begin{abstract}
A triple vector bundle is a cube of vector bundle structures which
commute in the (strict) categorical sense. A grid in a triple vector
bundle is a collection of sections of each bundle structure with
certain linearity properties. A grid provides two routes around 
each face of the triple vector bundle, and six routes from the base 
manifold to the total manifold; the warps measure the lack of 
commutativity of these routes. 

In this paper we first prove that the sum of the warps in a triple
vector bundle is zero. The proof we give is intrinsic and, we
believe, clearer than the proof using decompositions given earlier
by one of us. We apply this result to the triple tangent bundle
$T^3M$ of a manifold and deduce (as earlier) the Jacobi identity. 

We further apply the result to the triple vector bundle $T^2A$ 
for a vector bundle $A$ using a connection in $A$ to define a
grid in $T^2A$. In this case the curvature emerges from the
warp theorem. 
\end{abstract}

\newpage

\section{Introduction}
\label{sect:intro}

\subsection{Double vector bundles, grids and warps}

Double vector bundles arise naturally in Poisson geometry, 
in the connection theory of vector bundles, and generally
in the study of geometric objects with two compatible structures. 
Double vector bundles have been {\lq\lq}floating around'' 
since at least Dieudonn\'e's \cite{Dieudonne:III} treatment of 
connection theory, but the first systematic and general treatment  
was provided by Pradines \cite{Pradines77}. A recent account
with references is \cite[Chap.~9]{Mackenzie:GT}. We briefly
recall the necessary facts. 

A double vector bundle is a square of vector bundles as shown in 
the first figure of (\ref{eq:dvbs}). There are two vector bundle
structures on $D$, with bases $A$ and $B$, each of which is itself
a vector bundle on base $M$; the two structures on $D$ commute in
the categorical sense (see below), and the map $D\to A\times_MB$ formed by
the two bundle projections is a surjective submersion.
 \begin{equation}
 \label{eq:dvbs}
 \begin{tikzcd}[row sep=2cm, column sep = 2cm]
 D \arrow[r, "q^D_B"] \arrow[d,"q^D_A"] & B \arrow[d,"q_B"] \\
 A \arrow[r,"q_A"] &M,
 \end{tikzcd}
 \hspace*{12mm}
 \begin{tikzcd}[row sep=2cm, column sep = 2cm]
 TA \arrow[r, "T(q)"] \arrow[d,"p_A"] &TM \arrow[d,"p"] \\
 A\arrow[r,"q"] &M,
 \end{tikzcd}
 \hspace*{12mm}
 \begin{tikzcd}[row sep=2cm, column sep = 2cm]
 T^*A \arrow[r, ""] \arrow[d,"c_A"] &A^* \arrow[d,"q_*"] \\
 A\arrow[r,"q"] &M.
 \end{tikzcd}
 \end{equation}
The second and third figures in (\ref{eq:dvbs}) show two standard examples
arising from an arbitrary vector bundle $A$. If $A$ has a Poisson structure
then it is linear if and only if the associated map $T^*A\to TA$ is a morphism 
of double vector bundles (in an obvious sense) of the structures above. 
When this is so, $A$ is the dual of a Lie algebroid \cite{Courant:1990t}. 
The third structure was introduced in global form in \cite{MackenzieX:1994}.
We give more details on this \dvb in Example \ref{ex:MX}.
Double vector bundles also arise in the Lie theory of double Lie groupoids
\cite{Mackenzie:2011}; we will not consider this theory here.

Each element $d$ of a \dvb $D$ may be represented in outline by the diagram
in (\ref{eq:outlines}) which shows the projections of $d$ under the two
bundle projections. Given another element $d'$ as shown, the sum over $A$
has the outline shown in the third figure.  

\begin{equation}
\label{eq:outlines}
\begin{tikzcd}[row sep=1cm, column sep = 1cm]
d \arrow[r,mapsto," "]\arrow[d,mapsto, swap," "] 
&b \arrow[d,mapsto," "]\\
a\arrow[r,mapsto," "]
&m, 
\end{tikzcd}
\hspace*{15mm}
\begin{tikzcd}[row sep=1cm, column sep = 1cm]
d' \arrow[r,mapsto," "]\arrow[d,mapsto, swap," "] 
&b' \arrow[d,mapsto," "]\\
a\arrow[r,mapsto," "]
&m, 
\end{tikzcd}
\hspace*{15mm}
\begin{tikzcd}[row sep=7mm, column sep = 7mm]
d \plus{A} d'\arrow[r,mapsto," "]\arrow[d,mapsto, swap," "] 
&b+b' \arrow[d,mapsto," "]\\
a\arrow[r,mapsto," "]
&m. 
\end{tikzcd}
\end{equation}

The statement that the two vector bundle structures on $D$ `commute in
the categorical sense' implies for the additions that
\begin{equation}
\label{eq:ich}
(d_1\plus{A} d_2) \plus{B} (d_3\plus{A} d_4) =   
(d_1\plus{B} d_3) \plus{A} (d_2\plus{B} d_4), 
\end{equation}
where $(d_i; a_i, b_i; m)$, $i = 1, \dots, 4$, have 
$a_1 = a_2$, $a_3 = a_4$, $b_1 = b_3$ and $b_2 = b_4$. There are similar
conditions involving the scalar multiplications. 

It follows that for elements $d$ which project to zeros under both bundle
projections, the two additions, and the scalar multiplications, coincide. 
Under these operations the set of such elements forms a vector bundle over $M$, 
called the \emph{core of $D$} \cite{Pradines77}, usually denoted $C$. 

The core $C$ is a submanifold of $D$; every element of $C$ is an element of $D$.
When working with examples, the core can usually be identified with a familiar
vector bundle and it can be important to distinguish between elements of
this bundle and the corresponding element of the \dvb. For example, the core
of the \dvb $TA$, the middle diagram in (\ref{eq:dvbs}), can be identified
with $A$ itself, \cite[9.1.7]{Mackenzie:GT}. Therefore, an element $a$ of the
core $A$, can be viewed either as an element of $A$, or as an element of $TA$.
In the latter case, we denote it by $\overline{a}\in TA$. For general \dvbs and
\tvbs, this distinction is usually not necessary, so in sections \ref{sect:gaw},
\ref{sect:htseoE}, and \ref{sect:proof-of-theorem}, we will not write bars over
core elements. This distinction will be made clearly in section \ref{sect:T2E}. 

Now suppose that $(d;a,b;m)$ and $(d';a',b';m)$ have $a = a'$ and $b = b'$. 
Then there is a unique $c\in C$ such that 
\begin{equation}
\label{eq:ced}
d = d'\addition\limits_A(c\addition\limits_B
\tilde{0}^A_a) = 
d'\addition\limits_B
(c\addition\limits_A\tilde{0}^B_b).
\end{equation}
In equations of this type, what is usually important is that $d-d'$, 
calculated in either structure, gives the core element $c$ plus an
appropriate zero. We will indicate this by $d - d' \ced{}{c}$. 
We use this notation from subsection \ref{subsect:fwt} onwards. 

We now describe the original motivating example for the concepts of
grid and warp. 

In the 1988 edition of their book \cite[p.297]{AbrahamMR:MTAA},
Abraham, Marsden and Ra{\c t}iu gave the following formula for the Lie
bracket of vector fields $X$ and $Y$ on a manifold $M$, 
\begin{equation}
\label{eq:AMRi}
T(Y)(X(m)) - \wtilde{X}(Y(m)) = ([X,Y](m))^\uparrow(Y(m)),
\end{equation}
where $\wtilde{X}$ is the complete lift of $X$ to a vector field on $TM$ and the
uparrow denotes the vertical lift to $TM$ of the vector $[X,Y](m)$ to $Y(m)$. 
The complete lift, or tangent lift, 
$\wtilde{X}$ is $J\circ T(X)$ where $J\co T^2M\to T^2M$ is the canonical 
involution which interchanges the two bundle structures on $T^2M$.
The \dvb $T^2M$ is a special case of the middle diagram of 
(\ref{eq:dvbs}), where $A = TM$, and its core \vb is yet a third copy of $TM$.
The left hand side of (\ref{eq:AMRi}) is encapsulated in (\ref{eq:AMRii}).

\begin{equation}
\label{eq:AMRii}
\begin{tikzcd}[row sep=2cm, column sep = 2cm]
T^2M \arrow[r, "T(p)",swap]\arrow[r,<-, bend left = 25,"T(Y)"]
\arrow[d,"p_{TM}"]\arrow[d,<-, bend right = 25,swap, "\widetilde{X}"] &TM
\arrow[d,"p",swap]\arrow[d,<-, bend left = 25, "X"]\\
TM\arrow[r,"p"]\arrow[r,<-, bend right = 25,swap, "Y"]
&M.
\end{tikzcd} 
\end{equation}
If we look at the elements $T(Y)(X(m))$ and $\widetilde{X}(Y(m))$, we see that
they have the same outlines

\begin{equation*}
\begin{tikzcd}[row sep=1.5cm, column sep = 1.5cm]
T(Y)(X(m))\arrow[r,mapsto," "]\arrow[d,mapsto, swap," "] 
&X(m) \arrow[d,mapsto," "]\\
Y(m)\arrow[r,mapsto,swap," "]
&m, 
\end{tikzcd}\quad
\begin{tikzcd}[row sep=1.5cm, column sep = 1.5cm]
\widetilde{X}(Y(m))\arrow[r,mapsto," "]\arrow[d,mapsto, swap," "] 
&X(m) \arrow[d,mapsto," "]\\
Y(m)\arrow[r,mapsto,swap," "]
&m.
\end{tikzcd}
\end{equation*}

The two elements therefore determine a core element $c\in TM$. Taking $d = T(Y)(X(m))$
and $d' = \wtilde{X}(Y(m))$ in (\ref{eq:ced}), we have 
$$
T(Y)(X(m)) - \widetilde{X}(Y(m)) = \Bar{c} \addition\limits_{T(p)}
\widetilde{0}_{Y(m)}, $$
where the subtraction on the left is the usual subtraction of vectors which are 
tangent to $TM$ at $Y(m)$, and the addition on the right is addition in $T(p)\co
T^2M\to TM$.
That is, $\Bar{c} \addition\limits_{T(p)} \widetilde{0}_{Y(m)}$ is the vertical
lift of $c$ to $Y(m)$ and so, by (\ref{eq:AMRi}), $c = [X,Y](m)$. 

A comment on the notation of the last equation. In the case of a general \dvb
$D$, the two additions $\addition\limits_A$ and $\addition\limits_B$ are
distinct. In the case of $T^2M$ however, both side bundles are copies of $TM$.
To distinguish between the two additions, we use the projection maps, for
example, addition in $T^2M\xrightarrow{T(p)}TM$ will be denoted by
$\addition\limits_{T(p)}$. We adopt this notation whenever necessary, especially
in sections \ref{sect:T2E} and \ref{sect:Jac}.

Note that (\ref{eq:AMRi}) needs to be proved in local coordinates, or in terms 
of the action of vector fields on functions. The use of (\ref{eq:AMRii}) 
expresses the result in a compact conceptual way. 

We now express these results in the terms which will be used throughout 
the paper. Consider a \dvb $D$ as in (\ref{eq:dvbs}). 

\begin{definition}
A pair of sections $X\in\Ga A$ and $\xi\in \Ga_BD$ form a \emph{linear section}
of $D$ if $\xi$ is a morphism of vector bundles over $X$. 

A \emph{grid on $D$} is a pair of linear sections $(\xi, X)$ and $(\eta,Y)$
as shown in (\ref{eq:dgrid}).
\end{definition}

\begin{equation}
\label{eq:dgrid}
\begin{tikzcd}[row sep=1.5cm, column sep = 1.5cm]
D \arrow[r]\arrow[r,<-, bend left = 25, "\xi"]
\arrow[d]\arrow[d,<-,bend right = 25,swap, "\eta"] 
&B \arrow[d]\arrow[d,<-, bend left = 25, "Y"]\\
A\arrow[r]\arrow[r,<-, bend right = 25,swap, "X"]
&M.
\end{tikzcd}
\end{equation}
For each $m\in M$, $\xi(Y(m))$ and $\eta(X(m))$ have the same outline. 
They therefore determine an element of the core $C$ and, as $m$ varies, 
a section of $C$ which we denote $\warp(\xi, \eta)$. Precisely, 
\begin{equation}\label{eqn:warp of xi and eta}
\xi(Y(m))\subtraction\limits_A \eta(X(m)) = \warp(\xi,\eta)(m)\addition\limits_B
\tilde{0}_{X(m)},\qquad
\xi(Y(m))\subtraction\limits_B  \eta(X(m))= \warp(\xi,\eta)(m)\addition\limits_A
\tilde{0}_{Y(m)}.
\end{equation}

\begin{definition}
\label{df:warp}
The \emph{warp} of the grid consisting of $(\xi, X)$ and $(\eta,Y)$ is 
$\warp(\xi,\eta)\in\Ga C$. 
\end{definition}

Note that $\warp(\xi,\eta)$ changes sign if $\xi$ and $\eta$ are interchanged. 
Our convention gives the positive sign to the counterclockwise composition $\xi\circ Y$. 

The question of signs --- or orientations --- will haunt us throughout the paper. 
Later on, we will see that there are various rules 
that, in many cases, determine which difference to take as the positive warp. 
These rules generally follow from established conventions of differential geometry. 

Equation (\ref{eq:AMRi}) can now be expressed as saying that the warp of
(\ref{eq:AMRii}) is $[X,Y]$. 

Before proceeding, we give two further examples of grids and warps in \dvbs,
and an alternative formula. 

\begin{example}
\label{ex:conn}
Consider the \dvb $TA$, the middle diagram in (\ref{eq:dvbs}), where $(A,q,M)$ is
a \vb, and let $\nabla$ be a connection in $A$.

For a vector field $Z$ on $M$, denote by $Z^\hlift$ the horizontal lift
of $Z$ to $A$; the word `horizontal' here has its standard meaning in
connection theory, and does not refer to the structures in $TA$.
Then $Z^\hlift$ is a linear vector field over $Z$. 

Take any $\mu\in\Ga A$ and form the grid shown in (\ref{eq:GT}). 
Then the warp of the grid is $\nabla_Z\mu$; that is, for $m\in M$, 
\begin{equation}
\label{eq:conn}
T(\mu)(Z(m)) - Z^\hlift(\mu(m)) = ((\nabla_Z\mu)(m))^\uparrow(\mu(m)),
\end{equation}
where the right hand side is the vertical lift of $(\nabla_Z\mu)(m)\in T_mM$ to 
$T_{\mu(m)}A$. 

\begin{equation}
\label{eq:GT}
\begin{tikzcd}[row sep=2cm, column sep = 2cm]
TA \arrow[r,"T(q)",swap]\arrow[r,<-, bend left = 25,
"T(\mu)"]\arrow[d,"p_A"]\arrow[d,<-, bend right = 25,swap, "Z^{\hlift}"] &TM
\arrow[d,"p",swap]\arrow[d,<-, bend left = 25, "Z"]\\
A\arrow[r,"q"]\arrow[r,<-, bend right = 25,swap, "\mu"]
&M.
\end{tikzcd}
\end{equation}

For details see \cite[\S3.4]{Mackenzie:GT}. This example is central to \S\ref{sect:T2E}. 

Some historical remarks may be in order. 

The most usual global language for working with connections in vector bundles
is that of covariant derivatives, nowadays with the $\nabla$ notation. This 
formulation goes back to early work on surfaces and 
Riemannian geometry. It was formalized as a general abstract concept by 
Koszul in lectures given in 1960 \cite{Koszul:1960}. 

For connections in principal bundles a similarly convenient formulation 
is not available. Kobayashi and Nomizu \cite{KN1} gave two global definitions
of a connection in a principal bundle $P(M,G)$: as a suitable $\gog$-valued
$1$-form on $P$ and as an invariant horizontal distribution on $P$. The latter
defines, and is equivalent to, a lifting of vector fields on $M$ to invariant
horizontal vector fields on $P$. 

In the 1970s this lifting formulation was applied to connections in vector bundles.  
In modern language this approach defines a connection in $(A,q,M)$ to be a map 
$\X(M)\to\X(A)$, $Z\mapsto Z^\hlift$, which preserves addition, which sends $fZ$, 
for $f\in C^\infty(M)$, to $(f\circ q)Z^\hlift$, and which is such that 
$Z^\hlift$ is a linear vector field over $Z$. 

Dieudonn\'e \cite[XVII.16]{Dieudonne:III} 
and others expressed such a lifting in terms of a suitably linear
right-inverse to the map $TA\to TM\times_M A$ which combines the two
projections. Such a map also defines a corresponding left-inverse 
$TA\to A$ into the core, and this is the formulation which
Besse \cite{Besse:MAWGC} used. Equation (\ref{eq:conn}) may be
discerned on page~38 of \cite{Besse:MAWGC} and is a special case
of (17.17.2.1) in \cite{Dieudonne:III}. 

The important fact is that the usual definition of a connection in a
vector bundle as a covariant derivative is equivalent to a suitable lifting of 
vector fields from the base to the total space. We will comment on curvature
in \S\ref{sect:T2E}. 
\end{example}

\begin{example}
\label{ex:MX}
  Consider the \dvb $T^*A$, the third diagram in (\ref{eq:dvbs}), where $(A,q,M)$
  is a \vb.

  Given a section $\ph\in\Ga A^*$, denote by $\ell_\ph$ the corresponding linear
  function $A\to\R$. Then the $1$-form $d\ell_\ph$ is a linear section over $\ph$.

  Likewise given $\mu\in\Ga A$, we obtain a $1$-form $d\ell_\mu$ on $A^*$. Composing
  with the canonical diffeomorphism $R\co T^*(A^*)\to T^*A$ \cite{MackenzieX:1994},
  \cite[\S9.5]{Mackenzie:GT}, we obtain a linear section of $T^*A\to A^*$ over $\mu$. 

  It was proved in \cite{MackenzieX:1994} that
  $$
R(d\ell_\mu(\ph(m))) - d\ell_\ph(\mu(m)) = - q^*(d\langle\ph,\mu\rangle)(\mu(m)).
$$
(or see \cite[9.5.3]{Mackenzie:GT}). This shows that the warp of the grid consisting of
$(d\ell_\ph,\ph)$ and $(R\circ(d\ell_\mu),\mu)$ is $-d\langle\ph,\mu\rangle$.   
\end{example}

An alternative formula for the warp of a grid is given in an appendix, \S\ref{sect:app}. 

\subsection{Outline of the paper}

In a previous paper \cite{Mackenzie:pJihw} one of us used 
a grid in the triple vector bundle $T^3M$ to express the Jacobi identity 
as a statement about the warps of the grids in the constituent \dvbs. 
The proof given in that paper relied on a decomposition of $T^3M$ into seven 
copies of $TM$, and it was not clear whether the apparatus of grids and warps 
had provided a proof of the Jacobi identity or merely a formulation of it. 
One purpose of the present paper is to give an intrinsic proof of a general 
result for \tvbs and to resolve this question. 

Section \ref{sect:gaw} reviews the basic setup and notation for \tvbs. 
In subsection \ref{subsect:fwt} we formulate the main theorem of the
paper on the warps of a grid in a \tvb. 

In a \dvb two elements with the same outline determine a core element
(up to sign). 
In a \tvb there are intermediate levels at which two elements may have
the same outline, and there is more than one notion of core. Section 
\ref{sect:htseoE} is concerned with describing the core elements 
determined by pairs of elements for which some levels of the outlines
are equal. 

Section \ref{sect:proof-of-theorem} gives the proof of the warp theorem, 
Theorem~\ref{thm:grid}. This states, roughly, that given a grid in a
triple vector bundle $E$, the sum of the ultrawarps is zero. The ultrawarps
are the warps of the grids induced in the core double vector bundles by
the grid on $E$. The proof is intrinsic and does not rely on a decomposition
of $E$. 

In section \ref{sect:T2E} we consider a vector bundle $A$ and the triple
vector bundle $T^2A$. A connection in $A$ induces a grid in $T^2A$ and we
show that the concept of curvature arises from the warp theorem. 

Section~\ref{sect:Jac} presents the example of $T^3M$ and the deduction of
the Jacobi identity from the warp theorem. 

In an appendix, section~\ref{sect:app}, we show that the duality properties
of \dvbs give an alternative formula for the warp of a grid. 

In a future paper we will consider in detail the application of the
warp theorem to triple vector bundles with compatible bracket structures. 

\subsection*{Acknowledgements}

We are very grateful to Yvette Kosmann-Schwarzbach, Ted Voronov and Ping Xu
for valuable comments at various stages in the preparation of this paper.
We would also like to thank Fani Petalidou for thoroughly
reading an earlier version of the paper and catching many slips.

\newpage

\section{\Tvbs and the warp theorem}
\label{sect:gaw}

In this section we do two things. First, we set up
everything we need from \tvbs, in order to formulate the warp theorem 
(Theorem \ref{thm:grid}). 
Secondly, we describe the original formulation \cite{Mackenzie:pJihw}
of the theorem, and outline the steps which lead to an intrinsic proof.

\subsection{Basics on \tvbs}

We review the basic structure of \tvbs from 
\cite{Mackenzie:2005dts,Gracia-SazM:2009,Mackenzie:pJihw}. 

\begin{definition}
\label{df:tvb}
A \emph{\tvb} is a cube of vector bundle structures, as in (\ref{eq:tvb}), 
such that each face is a \dvb, such that the vector bundle operations 
in the upper faces are morphisms of \dvbs,  
and such that the decomposition condition described on 
page~\pageref{page:decomposition} below is satisfied. 
\end{definition}

\begin{equation}
\label{eq:tvb}
\begin{tikzcd}[row sep=1.5em,column sep=1.15em]
E_{1,2,3} \arrow[rr] \arrow[dr,swap,""] \arrow[dd,swap,""] &&
  E_{1,3} \arrow[dd,swap,"" near start] \arrow[dr,""] \\
& E_{2,3} \arrow[rr,crossing over,"" near start] &&
  E_3 \arrow[dd,""] \\
E_{1,2} \arrow[rr,"" near end] \arrow[dr,swap,""] && E_1
\arrow[dr,swap,""]
\\
& E_2 \arrow[rr,""] \arrow[uu,<-,crossing over,"" near end]&& M.
\end{tikzcd}
\end{equation}

In the rest of this section we work with a single \tvb $E$. 

By an \emph{upper face} we mean a face which has $E_{1,2,3}$
as total space. The \emph{lower faces} are the three faces which 
have $M$ as base manifold. We refer to the faces by the names
$$
\text{Back,\quad Front,\quad Left,\quad Right,\quad Up,\quad Down.}
$$
The total space of a triple vector bundle should be denoted, for 
consistency with the labelling scheme, by $E_{1,2,3}$ but we will 
usually denote it by $E$. 

How do we add elements in $E$? If $e,f\in E$ lie over the same point of $E_{2,3}$,
their sum has the outline shown in (\ref{eq:addinE}). 
\begin{equation}
\label{eq:addinE}
\begin{tikzcd}[row sep=1.15em,column sep=1.15em]
e \arrow[rr,""] \arrow[dr,swap,""] \arrow[dd,swap,""] &&
e_{1,3} \arrow[dd,swap,"" near start] \arrow[dr,""] \\
& e_{2,3} \arrow[rr,crossing over,"" near start] &&
e_3\arrow[dd,""] \\
e_{1,2} \arrow[rr,"" near end] \arrow[dr,swap,""] && e_1 \arrow[dr,swap,""]
\\
& e_2 \arrow[rr,""] \arrow[uu,<-,crossing over,"" near end]&& m
\end{tikzcd}\addition\limits_{2,3}
\begin{tikzcd}[row sep=1.15em,column sep=1.15em]
f \arrow[rr,""] \arrow[dr,swap,""] \arrow[dd,swap,""] &&
f_{1,3} \arrow[dd,swap,"" near start] \arrow[dr,""] \\
& e_{2,3} \arrow[rr,crossing over,"" near start] &&
e_3\arrow[dd,""] \\
f_{1,2} \arrow[rr,"" near end] \arrow[dr,swap,""] && f_1 \arrow[dr,swap,""]
\\
& e_2 \arrow[rr,""] \arrow[uu,<-,crossing over,"" near end]&& m
\end{tikzcd}=
\begin{tikzcd}[row sep=0.8em,column sep=1.15em]    
e\addition\limits_{2,3}f \arrow[rr,""] \arrow[dr,swap,""] \arrow[dd,swap,""] &&
e_{1,3}\addition\limits_{E_3} 
f_{1,3} \arrow[dd,swap,"" near start] \arrow[dr,""]
\\
& e_{2,3} \arrow[rr,crossing over,"" near start] &&
e_3\arrow[dd,""] \\
e_{1,2}\addition\limits_{E_2} 
f_{1,2} \arrow[rr,"" near end] \arrow[dr,swap,""]
&& e_1+f_1
\arrow[dr,swap,""]
\\
& e_2 \arrow[rr,""] \arrow[uu,<-,crossing over,"" near end]&& m.
\end{tikzcd}
\end{equation}
The outlines for scalar multiplication are similar.

\subsubsection*{Core \dvbs and the ultracore}

Since each face of $E$ is a \dvb, each face has a core \vb.

The cores of the lower faces $E_{i,j}$ are denoted $E_{ij}$ with the
comma removed. The core of the upper face with base manifold $E_k$ is
denoted $E_{ij,k}$. (This convention comes from \cite{Gracia-SazM:nfold-x}.)

Focus on the core \vbs of the Up and of the Down faces. The Up face 
projects to the Down face via the \dvb morphism which consists of the 
bundle projections $E_{1,2,3}\to E_{1,2}$, $E_{2,3}\to E_2$, $E_{1,3}\to E_1$ 
and $E_3\to M$. The restriction of $E_{1,2,3}\to E_{1,2}$ to $E_{12,3}$
goes into $E_{12}$ and inherits the vector bundle structure of $E_{1,2,3}\to E_{1,2}$. 
Together with the vector bundle structures on the cores of the Up face
and the Down face, this yields another \dvb, with total space $E_{12,3}$, 
which we call the \emph{(U-D) core \dvb}. 

Of course this can also be done for the other two pairs of parallel faces. 
So there are three core \dvbs, shown in (\ref{eq:UP-DOWN CORE DVB}).

\begin{equation}
\label{eq:UP-DOWN CORE DVB}
\qquad
\begin{tikzcd}[row sep=1.5cm, column sep = 1.5cm]
E_{23,1} \arrow[r,""]\arrow[d, swap,""] 
&E_{23} \arrow[d,""]\\
E_1\arrow[r,swap,""]
&M, 
\end{tikzcd}
\qquad
\begin{tikzcd}[row sep=1.5cm, column sep = 1.5cm]
E_{13,2} \arrow[r,""]\arrow[d, swap,""] 
&E_{13} \arrow[d,""]\\
E_2\arrow[r,swap,""]
&M,
\end{tikzcd}
\qquad
\begin{tikzcd}[row sep=1.5cm, column sep = 1.5cm]
E_{12,3} \arrow[r]
\arrow[d]&E_{12} 
\arrow[d]\\
E_3\arrow[r]
&M.
\end{tikzcd}
\end{equation}

Elements of the core of $E_{12,3}$ project to zeros in the Down face. 
In the Up face they project to zeros over the zero in $E_3$. It follows
that an element of the core of $E_{12,3}$ projects to zero in every bundle
structure. Equally the cores of the (B-F) and (L-R) \dvbs consist of 
the elements of $E_{1,2,3}$ which project to zeros in every bundle 
structure. Thus each \dvb in (\ref{eq:UP-DOWN CORE DVB}) has 
the same core. This is denoted $E_{123}$ (without commas) and called
the \emph{ultracore of $E$}.  

From the interchange laws it follows that the three 
additions on $E$, namely $\addition\limits_{1,2}$, $\addition\limits_{1,3}$, 
and $\addition\limits_{2,3}$, coincide on the ultracore and give it the
structure of a vector bundle over $M$. 

We can now state the decomposition condition needed in Definition~\ref{df:tvb}. 
\label{page:decomposition}

First consider a \dvb $D$ as in (\ref{eq:dvbs}), with core $C$. From $A$, $B$ 
and $C$ it is possible to define a \dvb structure on the pullback manifold 
$\Bar{D} = A\times_M B\times_M C$ by defining $\Bar{D}\to A$ to be the
pullback of $B\oplus C$ across $q_A\co A\to M$ and $\Bar{D}\to B$ to be the
pullback of $A\oplus C$ across $q_B\co B\to M$. This \dvb has the same side
bundles and same core as $D$, and the condition that the map $D\to A\times_M B$
be a surjective submersion implies that there are isomorphisms of \dvbs from
$D$ to $\Bar{D}$ which are the identity on $A$, $B$ and $C$. 

Now given a \tvb $E$, it is possible to build a \tvb $\Bar{E}$ from $E_1$, 
$E_2$, $E_3$, $E_{12}$, $E_{23}$, $E_{13}$ and $E_{123}$ by taking pullbacks
of Whitney sums in a similar way. The decomposition condition referred to
in Definition~\ref{df:tvb} is that there exists an isomorphism of \tvbs from
$E$ to $\Bar{E}$ which is the identity on the seven bundles listed above. 
See \cite{Gracia-SazM:2009} for full details. 

This condition ensures that in a \tvb $E$ (nontrival) grids, as defined
below always exist. 

\subsubsection*{Notation for zero sections}

The zero section of $E_1$ is denoted by $0^{E_1}:M\rightarrow E_1$, 
$m\mapsto 0^{E_1}_m$, with similar notations for $E_2$ and $E_3$.

The zero section of $E_{1,2}\rightarrow E_1$ is denoted by $\wtilde{0}^{1,2}:
E_1\rightarrow E_{1,2}$, $e_1\mapsto \wtilde{0}^{1,2}_{e_1}$. The
double zero of $E_{1,2}$ is denoted by $\odot^{1,2}_m$, with similar 
notations for the other \vb structures.

The zero section of $E\rightarrow E_{1,2}$ is denoted by
$\hat{0}:E_{1,2}\rightarrow E$, $e_{1,2}\mapsto \hat{0}_{e_{1,2}}$. Note that
the subscripts of the element $e_{1,2}$ are enough to indicate that this
is the zero section of $E$ over $E_{1,2}$; there is no need for 
superscripts on $\hat{0}$.

Finally, the triple zero of $E$ is denoted by $\odot^3_m$. 
This is the zero of the ultracore \vb. 

\subsubsection*{Grids in triple vector bundles}

A grid in a \dvb constitutes two linear sections.
In a \tvb the concept of grid requires what we call linear double sections. 

\begin{definition}
\label{df:du}
A \emph{down-up linear double section of $E$} is a collection of
sections 
\begin{equation*}
Z_{1,2}: E_{1,2}\rightarrow E_{1,2,3},\quad
Z_1: E_1\rightarrow E_{1,3},\quad
Z_2: E_2\rightarrow E_{2,3},\quad
Z: M\rightarrow E_3,
\end{equation*}
which form a morphism of \dvbs from the Down face to the Up face. 
\end{definition}
The core morphism of $Z_{1,2}$ defines a \vb morphism from the core of the Down
face to the core of the Up face. We denote this by $Z_{12}: E_{12}\rightarrow E_{12,3}$. 
It is a linear section over $Z: M\rightarrow E_3$.

In a similar fashion we define \emph{right-left} and \emph{front-back linear
double sections} of $E$.

\begin{definition}
\label{df:grid}
A \emph{grid on $E$} is a set of three linear double sections, one in each
direction, as shown in (\ref{eq:++}).   
\end{definition}
\begin{equation}
\label{eq:++}
\begin{tikzcd}[row sep=1.5em]
E_{1,2,3} \arrow[rr,<-,"Y_{1,3}"] \arrow[dr,<-,swap,"X_{2,3}"]
\arrow[dd,<-,swap,"Z_{1,2}"] && E_{1,3} \arrow[dd,<-,swap,"Z_1" near start]
\arrow[dr,<-,"X_3"] \\
& E_{2,3} \arrow[rr,<-,crossing over,swap,"Y_3" near start] &&
  E_3 \arrow[dd,<-,"Z"] \\
E_{1,2} \arrow[rr,<-,swap,"Y_1" near start] \arrow[dr,<-,swap,"X_2"] && E_1
\arrow[dr,<-,"X"]
\\
& E_2 \arrow[rr,<-,swap,"Y"] \arrow[uu,crossing over,"Z_2" near end]&& M.
\end{tikzcd} 
\end{equation}

That nontrivial grids always exist follows from the decomposability
condition. 

We note the following equations for future reference. 
They follow from the fact that the double sections are morphisms of \dvbs. 

For $e_{1,2}$, $e'_{1,2}$ over the same point of $E_1$, 
  \begin{equation}\label{eqn:Z-1,2 l.d.s:over-E-1}
  Z_{1,2}(e_{1,2}\subtraction\limits_{E_1} e'_{1,2}) =
Z_{1,2}(e_{1,2})\subtraction\limits_{1,3} Z_{1,2}(e'_{1,2}).
  \end{equation}
For $e_{1,2}$, $e'_{1,2}$ over the same point of $E_2$, 
  \begin{equation}\label{eqn:Z-1,2 l.d.s:over-E-2}
  Z_{1,2}(e_{1,2}\subtraction\limits_{E_2} e'_{1,2}) =
Z_{1,2}(e_{1,2})\subtraction\limits_{2,3} Z_{1,2}(e'_{1,2}).
  \end{equation}

\subsection{Formulation of the warp theorem}
\label{subsect:fwt}

We now have everything we need in order to describe the original formulation of
the theorem, as given in~\cite{Mackenzie:pJihw}.

Start with a grid on $E$. Focus on the Up face of the \tvb. We see that
$(Y_{1,3},Y_3)$ and $(X_{2,3},X_3)$ define a grid on the Up face. Denote its
warp by $\warp_{\text{up}}$. This is a section of the core \vb of the Up face,
that is, $\warp_{\text{up}}:E_3\rightarrow E_{12,3}$. Similarly for the Down
face, its warp $\warp_{\text{down}}$ is a section of the core \vb of the Down
face, so $\warp_{\text{down}}:M\rightarrow E_{12}$. It follows 
that $(\warp_{\text{up}},\warp_{\text{down}})$ is a 
linear section of the (U-D) core \dvb. Recall that the core morphism $Z_{12}$
of the linear double section $Z_{1,2}$ defines another linear section of the (U-D)
core \dvb. Therefore, we have the following grid on $E_{12,3}$:
\begin{equation*}
\begin{tikzcd}[row sep=1.5cm, column sep = 1.5cm]
E_{12,3} \arrow[r]\arrow[r,<-, bend left = 25, "Z_{12}"]
\arrow[d]\arrow[d,<-,bend right = 25,swap, "\warp_{\text{up}}"] &E_{12} 
\arrow[d]\arrow[d,<-, bend left =25, "\warp_{\text{down}}"]\\
E_3\arrow[r]\arrow[r,<-, bend right = 25,swap, "Z"]
&M.
\end{tikzcd}
\end{equation*}
We call the warp of this grid the \emph{Up-Down ultrawarp} and denote it by
$\uwarp_\ud$. It is a section of the ultracore $E_{123}$. 

Of course we can also build corresponding grids on the other two core \dvbs. 
We therefore have three ultrawarps, as shown in (\ref{eq:ultrawarps}).

\begin{equation}
\label{eq:ultrawarps}
\begin{tikzcd}[row sep=1.5cm, column sep = 1.5cm]
E_{23,1} \arrow[r]\arrow[r,<-, bend left = 25, "X_{23}"]\arrow[d]\arrow[d,<-,
bend right = 25,swap, "\warp_{\text{back}}"] &E_{23} \arrow[d]\arrow[d,<-, bend
left = 25, "\warp_{\text{front}}"]\\
E_1\arrow[r]\arrow[r,<-, bend right = 25,swap, "X"]
&M,
\end{tikzcd}\qquad
\begin{tikzcd}[row sep=1.5cm, column sep = 1.5cm]
E_{13,2} \arrow[r]\arrow[r,<-, bend left = 25, "Y_{13}"]\arrow[d]\arrow[d,<-,
bend right = 25,swap, "\warp_{\text{left}}"] &E_{13} \arrow[d]\arrow[d,<-, bend
left = 25, "\warp_{\text{right}}"]\\
E_2\arrow[r]\arrow[r,<-, bend right = 25,swap, "Y"]
&M,
\end{tikzcd}\qquad
\begin{tikzcd}[row sep=1.5cm, column sep = 1.5cm]
E_{12,3} \arrow[r]\arrow[r,<-, bend left = 25, "Z_{12}"]\arrow[d]\arrow[d,<-,
bend right = 25,swap, "\warp_{\text{up}}"] &E_{12} \arrow[d]\arrow[d,<-, bend left
= 25, "\warp_{\text{down}}"]\\
E_3\arrow[r]\arrow[r,<-, bend right = 25,swap, "Z"]
&M.
\end{tikzcd}
\end{equation}
We take the ultrawarps with the orientations opposite to (\ref{eq:ultrawarps});
that is, using a rough notation, 
\begin{equation}
\label{eq:core-orientations}
\uwarp_\baf: = 
\warp_{\text{back}}\circ X - X_{23}\circ \warp_{\text{front}},\quad 
\uwarp_\lr: = 
\warp_{\text{left}}\circ Y - Y_{13}\circ \warp_{\text{right}},\quad 
\uwarp_\ud: = 
\warp_{\text{up}}\circ Z - Z_{12}\circ \warp_{\text{down}}.
\end{equation}

We can now state the main theorem about grids in \tvbs. 
\begin{theorem}[Warp Theorem]
\label{thm:grid}
Given a \tvb $E$ and a grid in $E$ as in {\rm (\ref{eq:++})},     
\begin{equation}
\label{eq:grid}
\uwarp_\baf + \uwarp_\lr + \uwarp_\ud = 0.
\end{equation}
\end{theorem}

To give an intrinsic proof, we need to describe the ultrawarps in an alternative way.

So far, the only equation we have seen that describes the warp of a
grid on a \dvb is (\ref{eqn:warp of xi and eta}). Focus on the ultrawarp
$\uwarp_\ud$. From the grid on the (U-D) core \dvb, for $m\in M$, 
by (\ref{eqn:warp of xi and eta}) we have that
\begin{equation}\label{eqn:in-question}
(\warp_{\text{up}}\circ Z)(m)\subtraction\limits_{2,3}(Z_{12}\circ
\warp_{\text{down}})(m) = \hat{0}_{Z(m)} \addition\limits_{1,2}
\uwarp_\ud(m).
\end{equation}

How can we express $(\warp_{\text{up}}\circ Z)(m)$ and 
$Z_{12}(\warp_{\text{down}}(m))$ in a more useful way? 
About $\warp_{\text{up}}$, for any $e_3\in E_3$, again from
(\ref{eqn:warp of xi and eta}) we have that
\begin{equation*}
Y_{1,3}(X_3(e_3)) \subtraction\limits_{1,3} X_{2,3}(Y_3(e_3)) =
\hat{0}_{X_3(e_3)}  \addition\limits_{2,3} \warp_{\text{up}}(e_3).
\end{equation*}
Putting $e_3 = Z(m)$, we have
\begin{equation}
\label{eq:lambda-3a}
Y_{1,3}(X_3(Z(m)))\subtraction\limits_{1,3}
X_{2,3}(Y_3(Z(m))) =
\hat{0}_{X_3(Z(m))} \addition\limits_{2,3} \warp_{\text{up}}(Z(m)).
\end{equation}

We introduce a more succint notation, for use in calculations. 
\begin{equation}
\label{eq:123456}
\begin{split}
\one = Z_{1,2}(Y_1(X(m))), \quad
\two = Y_{1,3}(Z_1(X(m))), \quad
\thr = X_{2,3}(Z_2(Y(m))),\\
\fou = Z_{1,2}(X_2(Y(m))), \quad 
\fiv = Y_{1,3}(X_3(Z(m))), \quad
\six = X_{2,3}(Y_3(Z(m))).
\end{split}
\end{equation}
Now (\ref{eq:lambda-3a}) becomes
\begin{equation}\label{equation-lambda-3}
\fiv\subtraction\limits_{1,3}\six =
\hat{0}_{e_{1,3}'}\addition\limits_{2,3}\lambda_3,
\end{equation}
where $e_{1,3}' = X_3(Z(m))$ and $\lambda_3 = \warp_{\text{up}}(Z(m))$.
Note the following. In the case of a \dvb, we can rewrite 
(\ref{eqn:warp of xi and eta}) as
\begin{equation*}
\warp(\xi,\eta)(m) = 
((\xi\circ Y)(m)\subtraction\limits_{A}(\eta\circ X)(m))\subtraction\limits_{B}0^D_{X(m)}.
\end{equation*}
In the case of a \tvb, this is also possible. If we tried a similar calculation
on (\ref{equation-lambda-3}), since 
$\hat{0}_{e_{1,3}'}\subtraction\limits_{2,3}\hat{0}_{e_{1,3}'} =
 \hat{0}_{e_3}$, we would have
\begin{equation*}
(\fiv\subtraction\limits_{1,3}\six)\subtraction\limits_{2,3}\hat{0}_{e_{1,3}'} = 
\hat{0}_{e_3}\addition\limits_{2,3}\lambda_3.
\end{equation*}
Since $\hat{0}_{e_3}$ plays the role of the double zero of the Up face, over
$e_3$, we have $\hat{0}_{e_3}\addition\limits_{2,3}\lambda_3 =
\lambda_3$. So in total, we can rewrite (\ref{equation-lambda-3}) as
\begin{equation}\label{rewritten-equation-lambda-3}
\lambda_3 =
(\fiv\subtraction\limits_{1,3}\six)\subtraction\limits_{2,3}\hat{0}_{e_{1,3}'}.
\end{equation}
About $(Z_{1,2}\circ \warp_{\text{down}})(m)$, first write 
$\warp_{\text{down}}(m)$ out using (\ref{eqn:warp of xi and eta}) as
\begin{equation*}
Y_1(X(m))\subtraction\limits_{E_1}X_2(Y(m)) =
\wtilde{0}^{1,2}_{X(m)}\addition\limits_{E_2}\warp_{\text{down}}(m).
\end{equation*}
Apply $Z_{1,2}$ to this, and using (\ref{eqn:Z-1,2 l.d.s:over-E-1}) and
(\ref{eqn:Z-1,2 l.d.s:over-E-2}), it follows that
\begin{equation*}
Z_{1,2}(Y_1(X(m)))\subtraction\limits_{1,3}Z_{1,2}(X_2(Y(m))) = 
\hat{0}_{Z_1(X(m))}\addition\limits_{2,3}Z_{12}(\warp_{\text{down}}(m))
\end{equation*}
Again, for reasons of economy of space, rewrite this as
\begin{equation*}
\one\subtraction\limits_{1,3}\fou  = \hat{0}_{e_{1,3}}\addition\limits_{2,3}k_3, 
\end{equation*}
where $e_{1,3} = Z_1(X(m))$ and $k_3 = Z_{12}(\warp_{\text{down}}(m))$. 
Alternatively, as we did for $\lambda_3$
\begin{equation}\label{rewritten-equation-k-3}
k_3 =(\one\subtraction\limits_{1,3}\fou)
\subtraction\limits_{2,3}\hat{0}_{e_{1,3}}
\end{equation}
Let us go back to (\ref{eqn:in-question}). We can rewrite this as
\begin{equation*}
\lambda_3\subtraction\limits_{2,3}k_3 =
\hat{0}_{e_3}\addition\limits_{1,2}\uwarp_\ud(m),
\end{equation*}
and using (\ref{rewritten-equation-lambda-3}) and
(\ref{rewritten-equation-k-3}), we have that
\begin{equation*}
\left((\fiv\subtraction\limits_{1,3}\six)\subtraction\limits_{2,3}\hat{0}_{e_{1,3}'}\right)
\subtraction\limits_{2,3} 
\left((\one\subtraction\limits_{1,3}\fou)
\subtraction\limits_{2,3}\hat{0}_{e_{1,3}}\right) = 
\hat{0}_{e_3}\addition\limits_{1,2}\uwarp_\ud(m)
\end{equation*}
or, more elegantly, using interchange laws,
\begin{multline}
\label{eq:uUDelegant}
(\fiv\subtraction\limits_{1,3}\six)\subtraction\limits_{2,3}(\one\subtraction\limits_{1,3}\fou)
= (\hat{0}_{e_{1,3}'}\addition\limits_{2,3}\lambda_3)\subtraction\limits_{2,3}
(\hat{0}_{e_{1,3}}\addition\limits_{2,3}k_3)\\ = 
(\hat{0}_{e_{1,3}'}\subtraction\limits_{2,3}\hat{0}_{e_{1,3}})\addition\limits_{2,3}
(\lambda_3\subtraction\limits_{2,3}k_3) = 
(\hat{0}_{e_{1,3}'}\subtraction\limits_{2,3}\hat{0}_{e_{1,3}})\addition\limits_{2,3}
(\hat{0}_{e_3}\addition\limits_{1,2}\uwarp_\ud(m)).
\end{multline}
In calculations it is generally preferable to use equations of the form
(\ref{equation-lambda-3}), and to avoid equations of the form
(\ref{rewritten-equation-lambda-3}). 

Therefore, in order to describe ultrawarps such as $\uwarp_\ud(m)$, 
we will use equations of the form (\ref{eq:uUDelegant}) and we
will often use the abbreviated notation 
\begin{equation*}
(\fiv-\six)-(\one-\fou) \ced{} \uwarp_\ud(m),
\end{equation*}
introduced after (\ref{eq:ced}). 

It is worth emphasizing that the above arguments rely on the fact that core
and ultracore elements are uniquely determined by equations such as 
(\ref{eqn:warp of xi and eta}). 

There are similar abbreviated equations for the other two ultrawarps. Altogether we have
\begin{subequations}
\label{eq:uXYZ}
\begin{eqnarray}
(\one-\two)-(\thr-\six)\ced{} \uwarp_\baf(m),\label{simple-u1}\\
(\thr-\fou)-(\fiv-\two)\ced{} \uwarp_\lr(m),\label{simple-u2}\\
(\fiv-\six)-(\one-\fou)\ced{} \uwarp_\ud(m),\label{simple-u3}
\end{eqnarray}
\end{subequations}
and from now on we will use a further shortening of the notation
\begin{equation*}
\uwarp_\baf(m) = u_1,\qquad 
\uwarp_\lr(m) = u_2, \qquad
\uwarp_\ud(m) = u_3. 
\end{equation*}

The main difficulty in proving (\ref{eq:grid}) is that we cannot simply 
add and subtract the expressions in (\ref{eq:uXYZ}), since the operations 
are over different \vb structures. The apparatus of the next section 
overcomes this difficulty. 

\subsubsection*{Orientation}

A further problem arises from the fact that the warp of a grid 
in a \dvb is only defined up to sign. We now need to consider 
how to choose these signs consistently for a grid in a \tvb. 
This is a question of fixing orientations. 

We choose to orient each upper face so that the positive term in the 
formula for the warp defines the outward normal by the right-hand rule. 
We then take the positive and negative terms in the opposite lower face
to match those in the upper face; that is, we orient the lower faces 
so that the positive term in the warp defines the inward normal.

Thus the orientation of the Up face determines the signs in the first
subtraction in (\ref{simple-u3}) below and the orientation of the Down face
determines the signs in the second subtraction. 

The ``middle subtractions'' in (\ref{eq:uXYZ}), that is, the orientations 
of the core \dvbs, is an independent choice, equivalent to the choice of
signs in (\ref{eq:core-orientations}). What matters here is consistency: 
if we took all three warps with the opposite signs, that would also be fine.

\section{Preliminaries for the proof of the warp theorem}
\label{sect:htseoE}

This section contains the main technical work needed for the proof of the
warp theorem. We first describe our approach. 

We want to find how the three ultrawarps are related, more specifically, for 
$m\in M$, we want to find a relation between three ultracore elements, 
$u_1, u_2, u_3$. The best way
to do this, where slightly easier calculations are involved, is to manipulate
the differences of the six elements, (\ref{simple-u1}), (\ref{simple-u2}), and
(\ref{simple-u3}). That is exactly what we do in this and the following section. 
In this section we obtain formulas expressing the difference of two elements of a \tvb
in terms of core elements and zeros. If two elements of a \tvb can be subtracted
then their outlines must have at least one face in common. Cases where the outlines 
have two or more faces in common arise repeatedly in the rest of the paper. 
Each of these cases needs individual treatment.

Looked at from another point of view, two elements of a \tvb which can be subtracted
may admit exactly one, or two, or all three, of the subtractions 
$\subtraction\limits_{1,2}$, $\subtraction\limits_{1,3}$, and 
$\subtraction\limits_{2,3}$. 

\subsection{First case: two elements that have the same outline}
\label{ssect:sameoutline}

Let $e$ and $e'$ have exactly the same outline
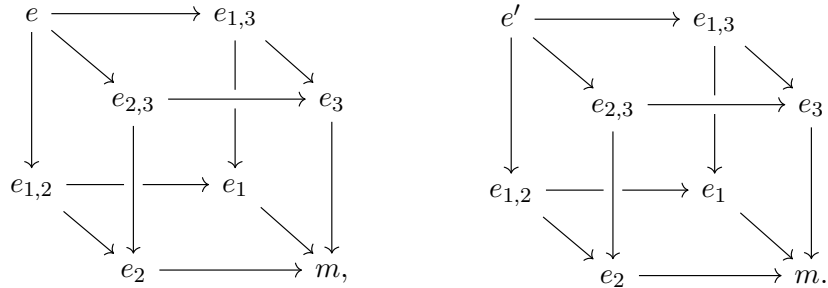
\begin{equation*}
\begin{tikzcd}[row sep=1.5em,column sep=1.15em]
e \arrow[rr] \arrow[dr,swap,""] \arrow[dd,swap,""] &&
  e_{1,3} \arrow[dd,swap,"" near start] \arrow[dr,""] \\
& e_{2,3} \arrow[rr,crossing over,"" near start] &&
  e_3 \arrow[dd,""] \\
e_{1,2} \arrow[rr,"" near end] \arrow[dr,swap,""] && e_1
\arrow[dr,swap,""]
\\
& e_2 \arrow[rr,""] \arrow[uu,<-,crossing over,"" near end]&& m,
\end{tikzcd}
\qquad\qquad
\begin{tikzcd}[row sep=1.5em,column sep=1.15em]
e' \arrow[rr] \arrow[dr,swap,""] \arrow[dd,swap,""] &&
  e_{1,3} \arrow[dd,swap,"" near start] \arrow[dr,""] \\
& e_{2,3} \arrow[rr,crossing over,"" near start] &&
  e_3 \arrow[dd,""] \\
e_{1,2} \arrow[rr,"" near end] \arrow[dr,swap,""] && e_1
\arrow[dr,swap,""]
\\
& e_2 \arrow[rr,""] \arrow[uu,<-,crossing over,"" near end]&& m.
\end{tikzcd}
\end{equation*}
All three differences $e\subtraction\limits_{1,2} e'$, 
$e\subtraction\limits_{1,3} e'$, $e\subtraction\limits_{2,3} e'$ 
are defined. 

\textbf{Step 1.}
Focus on the Back faces of $e$ and $e'$
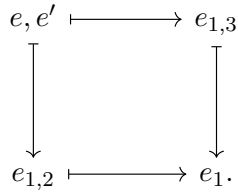
\begin{equation*}
\begin{tikzcd}[row sep=1.5cm, column sep = 1.5cm]
e,e' \arrow[r,mapsto,""]\arrow[d,mapsto, swap,""] 
&e_{1,3} \arrow[d,mapsto,""]\\
e_{1,2}\arrow[r,swap,mapsto,""]
&e_1.
\end{tikzcd}
\end{equation*}
Then, from \dvb theory, we can write
\begin{equation*}
e\subtraction\limits_{1,2} e' = k_1\addition\limits_{1,3}
\hat{0}_{e_{1,2}},\qquad
e\subtraction\limits_{1,3} e' = k_1\addition\limits_{1,2}
\hat{0}_{e_{1,3}},
\end{equation*}
where $k_1\in E_{23,1}$, the core of the Back face, with outline
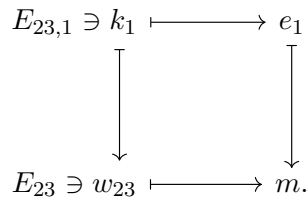
\begin{equation*}
\begin{tikzcd}[row sep=1.5cm, column sep = 1.5cm]
E_{23,1}\ni k_1 \arrow[r,mapsto,""]\arrow[d, swap,"",mapsto,xshift = 6mm] 
&e_1 \arrow[d,mapsto,""]\\
E_{23}\ni w_{23}\arrow[r,mapsto,swap,""]
&m.
\end{tikzcd}
\end{equation*}
\textbf{Step 2.} Show that $w_{23} = \odot^{2,3}_m$.

Use the morphism $q_{2,3}:E\rightarrow E_{2,3}$.
We know that 
$q_{2,3}(e\subtraction\limits_{1,3} e') = \tilde{0}^{2,3}_{e_3}$
and
\begin{equation*}
q_{2,3}(k_1\addition\limits_{1,2}\hat{0}_{e_{1,3}}) =
q_{2,3}(k_1)\addition\limits_{E_2}q_{2,3}(\hat{0}_{e_{1,3}}) =
w_{23}\addition\limits_{E_2}\tilde{0}^{2,3}_{e_3}. 
\end{equation*}
Therefore
\begin{equation*}
w_{23}\addition\limits_{E_2}\tilde{0}^{2,3}_{e_3} =
\tilde{0}^{2,3}_{e_3},
\end{equation*}
and, from \dvb theory, 
we have that $w_{23} = \tilde{0}^{2,3}_{0^{E_2}_m} = \odot^{2,3}_m$.
So $k_1$ has the outline
\begin{equation*}
\begin{tikzcd}[row sep=1.5cm, column sep = 1.5cm]
k_1 \arrow[r,mapsto,""]\arrow[d,mapsto, swap,""] 
&e_1 \arrow[d,mapsto,""]\\
\odot^{2,3}_m\arrow[r,swap,mapsto,""]
&m.
\end{tikzcd}
\end{equation*}
\textbf{Step 3.}
Applying \dvb theory again, we get 
\begin{equation*}
k_1 = u_1\addition\limits_{2,3}\hat{0}_{e_1},
\end{equation*}
where $u_1$ is an ultracore element.

\textbf{Step 4.}
Apply the same procedure to Left and Up faces of $e$ and $e'$. 

Focus on the Left faces of $e$ and $e'$
\begin{equation*}
  e\subtraction\limits_{2,3} e' =
  k_2\addition\limits_{1,2} \hat{0}_{e_{2,3}}, \qquad
  e\subtraction\limits_{1,2} e' =
  k_2\addition\limits_{2,3} \hat{0}_{e_{1,2}},
  \end{equation*}
where $k_2\in E_{13,2}$, core of Left face, with outline
\begin{equation*}
\begin{tikzcd}[row sep=1.5cm, column sep = 1.5cm]
k_2 \arrow[r,mapsto,""]\arrow[d, swap,mapsto,""] 
&e_2 \arrow[d,mapsto,""]\\
\odot^{1,3}_m\arrow[r,mapsto,swap,""]
&m.
\end{tikzcd}
\end{equation*}
So, we can write
\begin{equation*}
k_2 = u_2\addition\limits_{1,3} \hat{0}_{e_2},
\end{equation*}
where $u_2$ is an ultracore element.

Similarly for the Up faces
\begin{equation*}
  e\subtraction\limits_{1,3} e' =
  k_3\addition\limits_{2,3} \hat{0}_{e_{1,3}}, \qquad
  e\subtraction\limits_{2,3} e' =
  k_3\addition\limits_{1,3} \hat{0}_{e_{2,3}},
  \end{equation*}
where $k_3\in E_{12,3}$, core of Up face, with outline
\begin{equation*}
\begin{tikzcd}[row sep=1.5cm, column sep = 1.5cm]
k_3 \arrow[r,mapsto,""]\arrow[d, swap,mapsto,""] 
&e_3 \arrow[d,mapsto,""]\\
\odot^{1,2}_m\arrow[r,swap,mapsto,""]
&m,
\end{tikzcd}
\end{equation*}
so 
\begin{equation*}
k_3 = u_3\addition\limits_{1,2}\hat{0}_{e_3}
\end{equation*}
with $u_3$ an ultracore element.

\textbf{Step 5.} Show that $u_1 = u_2 = u_3$.

We show that $u_1 = u_3$. So far, we have two expressions for
$e\subtraction\limits_{1,3}e'$, namely:
\begin{equation}\label{eqn:equality-of-u1-and-u3}
k_1\addition\limits_{1,2}\hat{0}_{e_{1,3}} = 
k_3\addition\limits_{2,3}\hat{0}_{e_{1,3}}.
\end{equation}
Expand the left hand side of (\ref{eqn:equality-of-u1-and-u3}), mimicking the \dvb case:
\begin{equation*}
\hat{0}_{e_{1,3}} \addition\limits_{1,2}
(\hat{0}_{e_1}\addition\limits_{2,3}u_1) = 
(\hat{0}_{e_{1,3}} \addition\limits_{2,3}\hat{0}_{e_3})
 \addition\limits_{1,2}
(\hat{0}_{e_1}\addition\limits_{2,3}u_1) = 
(\hat{0}_{e_{1,3}} \addition\limits_{1,2}\hat{0}_{e_1})
 \addition\limits_{2,3}
(\hat{0}_{e_3}\addition\limits_{1,2}u_1)=
\hat{0}_{e_{1,3}}\addition\limits_{2,3}
(\hat{0}_{e_3}\addition\limits_{1,2}u_1).
\end{equation*}
Therefore, we see that (\ref{eqn:equality-of-u1-and-u3}) can be rewritten as:
\begin{equation*}
\hat{0}_{e_{1,3}}\addition\limits_{2,3}
(\hat{0}_{e_3}\addition\limits_{1,2}u_1) = 
\hat{0}_{e_{1,3}}\addition\limits_{2,3}
(\hat{0}_{e_3}\addition\limits_{1,2}u_3),
\end{equation*}
from where it follows that $u_1 = u_3$. Similarly, we can show that
$u_2 = u_3$.

At this point write $u_1= u_2=u_3$ to be $u$.

\textbf{Step 6.} We obtain six formulas for the differences between $e$ and $e'$. 

\begin{proposition}
\label{prop:sameoutline}
With the above notation, two elements $e$ and $e'$ which have the same outline
are related by   
\begin{eqnarray}
  e\subtraction\limits_{1,3} e'
  = \hat{0}_{e_{1,3}}\addition\limits_{1,2} (\hat{0}_{e_1}\addition\limits_{2,3}u)
  = \hat{0}_{e_{1,3}}\addition\limits_{2,3} (\hat{0}_{e_3}\addition\limits_{1,2}u), \nonumber\\
\label{eq:sameoutline}
e\subtraction\limits_{1,2} e'
= \hat{0}_{e_{1,2}}\addition\limits_{1,3} (\hat{0}_{e_1}\addition\limits_{2,3}u)
= \hat{0}_{e_{1,2}}\addition\limits_{2,3} (\hat{0}_{e_2}\addition\limits_{1,3}u),\\
e\subtraction\limits_{2,3} e'
= \hat{0}_{e_{2,3}}\addition\limits_{1,3} (\hat{0}_{e_3}\addition\limits_{1,2}u)
= \hat{0}_{e_{2,3}}\addition\limits_{1,2}
(\hat{0}_{e_2}\addition\limits_{1,3}u). \nonumber
\end{eqnarray}
\end{proposition}

What is important here is that the subtraction 
with respect to each structure results in the same ultracore element $u$. 

We will use the following special case later on. 

\begin{example}
\label{ex:13,12}
 If $e,e'$ are in one of the core \dvbs, for example if $e,e'\in E_{23,1}$, with
 outline
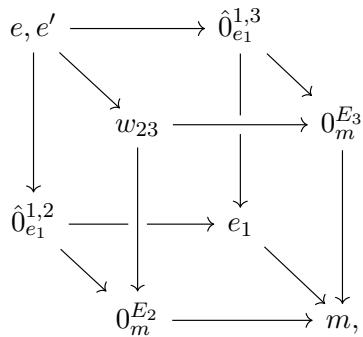
\begin{equation*}
\begin{tikzcd}[row sep=1.5em,column sep=1.15em]
e,e' \arrow[rr] \arrow[dr,swap,""] \arrow[dd,swap,""] &&
  \hat{0}_{e_1}^{1,3} \arrow[dd,swap,"" near start] \arrow[dr,""] \\
& w_{23} \arrow[rr,crossing over,"" near start] &&
  0^{E_3}_m \arrow[dd,""] \\
\hat{0}_{e_1}^{1,2} \arrow[rr,"" near end] \arrow[dr,swap,""] && e_1
\arrow[dr,swap,""]
\\
& 0^{E_2}_m \arrow[rr,""] \arrow[uu,<-,crossing over,"" near end]&& m,
\end{tikzcd}
\end{equation*}
then from {\rm (\ref{eq:sameoutline})} we have

\begin{equation*}
e\subtraction\limits_{2,3} e' = \hat{0}_{w_{23}}
\addition\limits_{1,3}
(\hat{0}_{0^{E_3}_m}\addition\limits_{1,2}u) 
= \hat{0}_{w_{23}}\addition\limits_{1,3}(\odot^3_m\addition\limits_{1,2}u)
= \hat{0}_{w_{23}}\addition\limits_{1,3}u
\end{equation*}
and
\begin{equation*}
e\subtraction\limits_{2,3} e' = \hat{0}_{w_{23}}
\addition\limits_{1,2} (\hat{0}_{0^{E_2}_m}\addition\limits_{1,3}u)=
\hat{0}_{w_{23}}\addition\limits_{1,2}(\odot^3_m\addition\limits_{1,3}u)
= \hat{0}_{w_{23}}\addition\limits_{1,2}u, 
\end{equation*}
and therefore
\begin{equation}
\label{eq:13,12}
 \hat{0}_{w_{23}}\addition\limits_{1,3}u = 
 \hat{0}_{w_{23}}\addition\limits_{1,2}u.
\end{equation}
\end{example}

\subsection{Second case: two elements that have two lower faces in common}

What happens if $e$ and $e'$ have only two of the lower faces in common?
Then only two of the three subtractions are defined. 
There are three cases to consider, each of which arises later. 

\subsubsection*{If $e$ and $e'$ have the same Front and Right face}
Since $e$ and $e'$ have the same $e_{1,3}$ and $e_{2,3}$, it follows
that they have the same $e_1$, $e_2$ and $e_3$.
However $e$ and $e'$ will differ at $e_{1,2}$ and $e_{1,2}'$, and these
will differ by a core element $w_{12}\in E_{12}$ of the core of the
Down face, that is
\begin{equation*}
  e_{1,2}\subtraction\limits_{E_1} e_{1,2}' =
  w_{12}\addition\limits_{E_2}\tilde{0}^{1,2}_{e_1},\qquad
  e_{1,2}\subtraction\limits_{E_2} e_{1,2}' =
  w_{12}\addition\limits_{E_1}\tilde{0}^{1,2}_{e_2}.
  \end{equation*}
It is useful to write out the outlines of these differences
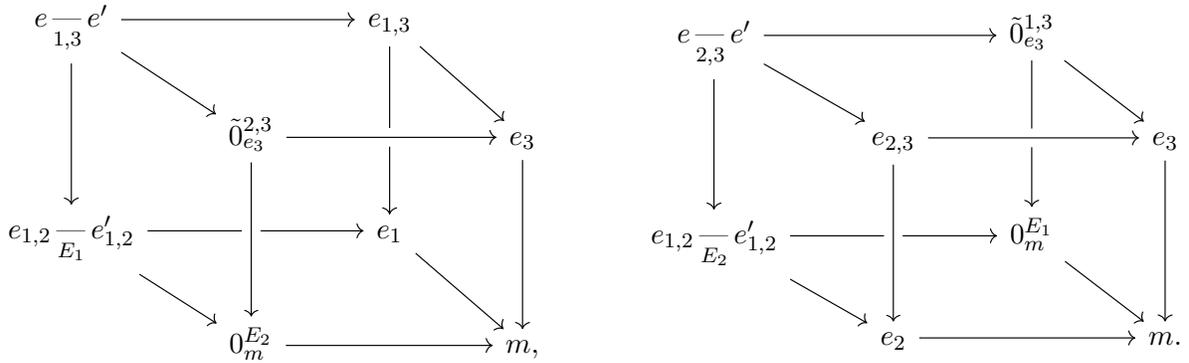
\begin{equation*}
\begin{tikzcd}[row sep=1.5em] 
e\subtraction\limits_{1,3} e' \arrow[rr]
\arrow[dr,swap,""] \arrow[dd,swap,""] && e_{1,3} \arrow[dd,swap,"" near start]
\arrow[dr,""] \\
& \tilde{0}^{2,3}_{e_3} \arrow[rr,crossing over,"" near start] &&
e_3 \arrow[dd,""] \\
e_{1,2}\subtraction\limits_{E_1} e_{1,2}' \arrow[rr,"" near end]
\arrow[dr,swap,""] && e_1
\arrow[dr,swap,""]
\\
& 0^{E_2}_m \arrow[rr,""] \arrow[uu,<-,crossing over,"" near end]&& m,
\end{tikzcd}
\qquad\quad
\begin{tikzcd}[row sep=1.5em] 
e\subtraction\limits_{2,3} e' \arrow[rr]
\arrow[dr,swap,""] \arrow[dd,swap,""] && \tilde{0}^{1,3}_{e_3}
\arrow[dd,swap,"" near start] \arrow[dr,""] \\
& e_{2,3}\arrow[rr,crossing over,"" near start] &&
e_3 \arrow[dd,""] \\
e_{1,2}\subtraction\limits_{E_2} e_{1,2}' \arrow[rr,"" near end]
\arrow[dr,swap,""] && 0^{E_1}_m
\arrow[dr,swap,""]
\\
& e_2 \arrow[rr,""] \arrow[uu,<-,crossing over,"" near end]&& m.
\end{tikzcd}
\end{equation*}
Since $e$ and $e'$ have the same Up face, again by applying \dvb theory, we can
write
\begin{equation}\label{same_front_and_right_faces_differences}
e\subtraction\limits_{1,3} e' = k\addition\limits_{2,3}
\hat{0}_{e_{1,3}},\qquad
 e\subtraction\limits_{2,3} e' =
k\addition\limits_{1,3} \hat{0}_{e_{2,3}},
\end{equation}
where $k\in E_{23,1}$, the core of the Up face.

Also, using the morphism $q_{1,2}: E\rightarrow E_{1,2}$, we show that
$q_{1,2}(k) = w_{12}$. First, 
\begin{equation*}
q_{1,2}(e\subtraction\limits_{1,3} e') =  e_{1,2}\subtraction\limits_{E_1} e_{1,2}' =
  w_{12}\addition\limits_{E_2}\tilde{0}^{1,2}_{e_1},
\end{equation*}
and 
\begin{equation*}
q_{1,2}(k\addition\limits_{2,3}\hat{0}_{e_{1,3}}) = 
q_{1,2}(k)\addition\limits_{E_2}\tilde{0}^{1,2}_{e_1},
\end{equation*}
hence $q_{1,2}(k) = w_{12}$, where $k$ has outline
\begin{equation*}
\begin{tikzcd}[row sep=1.5cm, column sep = 1.5cm]
E_{12,3}\ni k \arrow[r,mapsto,""]\arrow[d,mapsto, swap,"",xshift = 6mm] 
&e_3 \arrow[d,mapsto,""]\\
E_{12}\ni w_{12}\arrow[r,swap,mapsto,""]
&m.
\end{tikzcd}
\end{equation*}

\subsubsection*{If $e$ and $e'$ have the same Front and Down face}
In this case, the elements $e_{1,3}$ and $e_{1,3}'$ differ by a core
element $w_{13}$ of $E_{13}$
\begin{equation*}
e_{1,3}\subtraction\limits_{E_1} e_{1,3}' =
w_{13}\addition\limits_{E_3}\tilde{0}^{1,3}_{e_1},\qquad
e_{1,3}\subtraction\limits_{E_3} e_{1,3}' =
w_{13}\addition\limits_{E_1}\tilde{0}^{1,3}_{e_3}.
\end{equation*}
As before, we can write
\begin{equation}\label{same_front_and_bottom_faces_differences}
e\subtraction\limits_{1,2}e' = k\addition\limits_{2,3}
\hat{0}_{e_{1,2}},\qquad
e\subtraction\limits_{2,3}e' = k\addition\limits_{1,2}
\hat{0}_{e_{2,3}},
\end{equation}
with $k$ an element of the core of the Left face with outline
\begin{equation*}
\begin{tikzcd}[row sep=1.5cm, column sep = 1.5cm]
E_{13,2}\ni k \arrow[r,mapsto,""]\arrow[d, swap,mapsto,"",xshift = 6mm] 
&e_2 \arrow[d,mapsto,""]\\
E_{13}\ni w_{13}\arrow[r,swap,mapsto,""]
&m.
\end{tikzcd}
\end{equation*}
\subsubsection*{If $e$ and $e'$ have the same Right and Down face}
In this case, $e_{2,3}$ and $e_{2,3}'$ will differ by an element
$w_{23}\in E_{23}$ of the core of the Front face
\begin{equation*}
e_{2,3}\subtraction\limits_{E_2} e_{2,3}' =
w_{23}\addition\limits_{E_3}\tilde{0}^{2,3}_{e_2},\qquad
e_{2,3}\subtraction\limits_{E_3} e_{2,3}' =
w_{23}\addition\limits_{E_2}\tilde{0}^{2,3}_{e_3},
\end{equation*}
and as before
\begin{equation}\label{same_right_and_bottom_faces_differences}
e\subtraction\limits_{1,2}e' = k\addition\limits_{1,3}
\hat{0}_{e_{1,2}},\qquad
e\subtraction\limits_{1,3}e' = k\addition\limits_{1,2}
\hat{0}_{e_{1,3}},
\end{equation}
where $k$ is an element of the core of the Back face with outline
\begin{equation*}
\begin{tikzcd}[row sep=1.5cm, column sep = 1.5cm]
E_{23,1}\ni k \arrow[r,mapsto,""]\arrow[d, swap,"",mapsto,xshift = 6mm] 
&e_1\arrow[d,mapsto,""]\\
E_{23}\ni w_{23}\arrow[r,swap,mapsto,""]
&m.
\end{tikzcd}
\end{equation*}

\subsection{Special case: differences between zero elements}

Using the fact that $\hat{0}$ is the zero section
\begin{equation*}
\hat{0}_{e_{2,3}}\addition\limits_{1,3}\hat{0}_{e'_{2,3}} =
\hat{0}_{e_{2,3}\addition\limits_{E_3}e'_{2,3}},\qquad
\hat{0}_{e_{2,3}}\addition\limits_{1,2}\hat{0}_{e'_{2,3}} =
\hat{0}_{e_{2,3}\addition\limits_{E_2}e'_{2,3}}, 
\end{equation*}
and  $(-1)\mult\limits_{1,3}\hat{0}_{e_{2,3}} =
\hat{0}_{f_{2,3}}$ where $f_{2,3} = \subtraction\limits_{E_3}e_{2,3}$.

If we have two elements $e_{2,3}$ and $e_{2,3}'$ of $E_{2,3}$ that differ by a
core element $w_{23}\in E_{23}$, then
\begin{equation}
\label{eq:w23-1}
  e_{2,3}\subtraction\limits_{E_2}e'_{2,3} =
  w_{23}\addition\limits_{E_3} \tilde{0}^{2,3}_{e_2},\qquad
  e_{2,3}\subtraction\limits_{E_3}e'_{2,3} =
  w_{23}\addition\limits_{E_2} \tilde{0}^{2,3}_{e_3},
  \end{equation}
and the differences we are interested in are
  \begin{equation}
\label{eq:w23-2}
  \hat{0}_{e_{2,3}}\subtraction\limits_{1,2} \hat{0}_{e'_{2,3}} = 
  \hat{0}_{e_{2,3}}\addition\limits_{1,2}
  \hat{0}_{\subtraction\limits_{E_2}e'_{2,3}}=
  \hat{0}_{e_{2,3}\subtraction\limits_{E_2}e'_{2,3}} = 
  \hat{0}_{w_{23}\addition\limits_{E_3} \tilde{0}^{2,3}_{e_2}}=
  \hat{0}_{w_{23}}\addition\limits_{1,3}\hat{0}_{e_2},
  \end{equation}
  and
  \begin{equation}
\label{eq:w23-3}
  \hat{0}_{e_{2,3}}\subtraction\limits_{1,3} \hat{0}_{e'_{2,3}} = 
  \hat{0}_{e_{2,3}}\addition\limits_{1,3}
  \hat{0}_{\subtraction\limits_{E_3}e'_{2,3}}=
  \hat{0}_{e_{2,3}\subtraction\limits_{E_3}e'_{2,3}} = 
  \hat{0}_{w_{23}\addition\limits_{E_2} \tilde{0}^{2,3}_{e_3}}=
  \hat{0}_{w_{23}}\addition\limits_{1,2}\hat{0}_{e_3}.
  \end{equation}

There are similar formulas in the other two cases. 

\newpage

\section{Proof of the warp theorem}
\label{sect:proof-of-theorem}

\subsection{Notation}\label{subsect:notation-of-six-elements}

In this section we prove Theorem \ref{thm:grid}. We will use the
notation of (\ref{eq:123456}). We further simplify the notation
for elements of the lower faces and edges, as follows
\begin{gather*}
X(m):= e_1,\qquad Y(m):= e_2, \qquad Z(m):= e_3,\\
Z_1(X(m)): = e_{1,3},\quad  X_3(Z(m)): =e_{1,3}',\quad Z_2(Y(m)): = e_{2,3},\\ 
Y_3(Z(m)):=e_{2,3}',\quad Y_1(X(m)): = e_{1,2},\quad X_2(Y(m)):= e_{1,2}'.
\end{gather*}

The outlines of the elements in (\ref{eq:123456}) are now written as follows
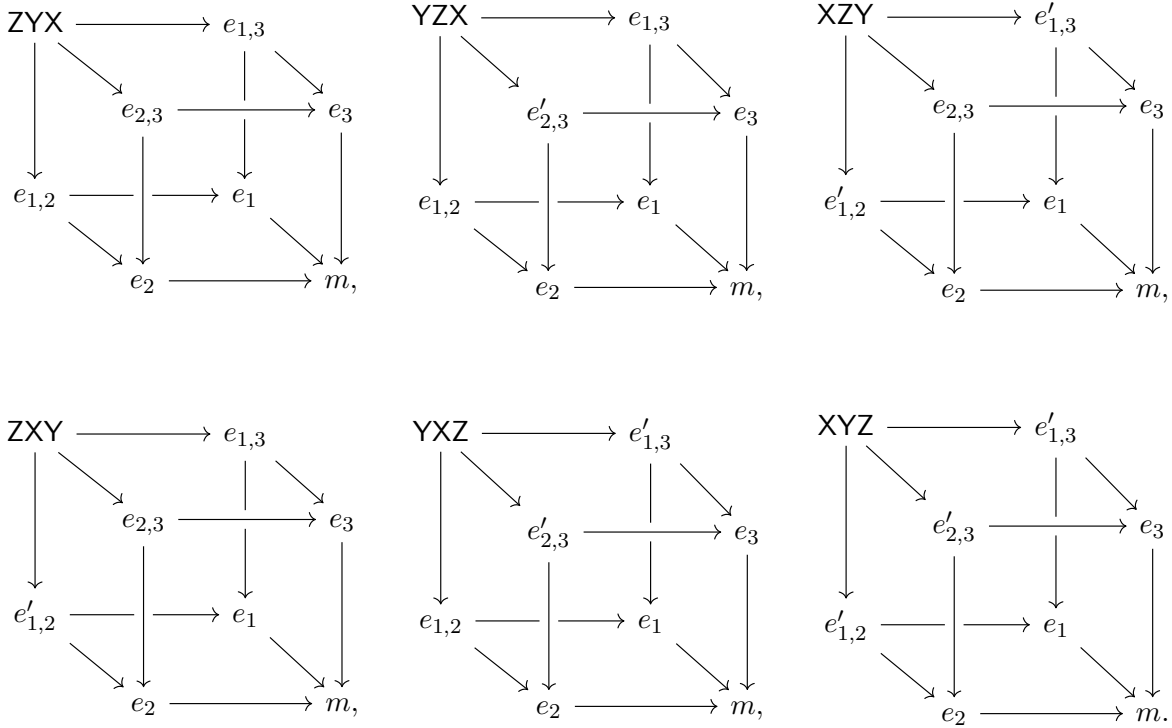
\begin{equation*}
\begin{tikzcd}[row sep=1.5em, column sep=1.15em]
\one \arrow[rr] \arrow[dr,swap,""]
\arrow[dd,swap,""] && e_{1,3} \arrow[dd,swap,"" near start] \arrow[dr,""] \\
& e_{2,3} \arrow[rr,crossing over,"" near start] &&
  e_3 \arrow[dd,""] \\
e_{1,2} \arrow[rr,"" near end] \arrow[dr,swap,""] && e_1
\arrow[dr,swap,""]
\\
& e_2 \arrow[rr,""] \arrow[uu,<-,crossing over,"" near end]&& m,
\end{tikzcd}
\quad
\begin{tikzcd}[row sep=1.5em, column sep=1.15em]
\two \arrow[rr] \arrow[dr,swap,""]
\arrow[dd,swap,""] && e_{1,3} \arrow[dd,swap,"" near start] \arrow[dr,""] \\
& e_{2,3}' \arrow[rr,crossing over,"" near start] &&
  e_3 \arrow[dd,""] \\
e_{1,2} \arrow[rr,"" near end] \arrow[dr,swap,""] && e_1
\arrow[dr,swap,""]
\\
& e_2 \arrow[rr,""] \arrow[uu,<-,crossing over,"" near end]&& m,
\end{tikzcd} 
\quad
\begin{tikzcd}[row sep=1.5em, column sep=1.15em]
\thr \arrow[rr] \arrow[dr,swap,""]
\arrow[dd,swap,""] && e_{1,3}' \arrow[dd,swap,"" near start]
\arrow[dr,""]\\
& e_{2,3} \arrow[rr,crossing over,"" near start] &&
 e_3\arrow[dd,""] \\
e_{1,2}' \arrow[rr,"" near end] \arrow[dr,swap,""] && e_1
\arrow[dr,swap,""]
\\
& e_2 \arrow[rr,""] \arrow[uu,<-,crossing over,"" near end]&& m,
\end{tikzcd}
\end{equation*}
\begin{equation*}
\begin{tikzcd}[row sep=1.5em, column sep=1.15em]
\fou \arrow[rr] \arrow[dr,swap,""]
\arrow[dd,swap,""] && e_{1,3} \arrow[dd,swap,"" near start]
\arrow[dr,""]
\\
& e_{2,3} \arrow[rr,crossing over,"" near start] &&
  e_3 \arrow[dd,""] \\
e_{1,2}' \arrow[rr,"" near end] \arrow[dr,swap,""] && e_1
\arrow[dr,swap,""]
\\
& e_2 \arrow[rr,""] \arrow[uu,<-,crossing over,"" near end]&& m,
\end{tikzcd}
\quad
\begin{tikzcd}[row sep=1.5em, column sep=1.15em]
\fiv \arrow[rr] \arrow[dr,swap,""]
\arrow[dd,swap,""] && e_{1,3}' \arrow[dd,swap,"" near start] \arrow[dr,""] \\
& e_{2,3}' \arrow[rr,crossing over,"" near start] &&
  e_3 \arrow[dd,""] \\
e_{1,2} \arrow[rr,"" near end] \arrow[dr,swap,""] && e_1
\arrow[dr,swap,""]
\\ 
& e_2\arrow[rr,""] \arrow[uu,<-,crossing over,"" near end]&& m,
\end{tikzcd}
\quad
\begin{tikzcd}[row sep=1.5em, column sep=1.15em]
\six \arrow[rr] \arrow[dr,swap,""]
\arrow[dd,swap,""] && e_{1,3}' \arrow[dd,swap,"" near start] \arrow[dr,""] \\
& e_{2,3}' \arrow[rr,crossing over,"" near start] &&
  e_3 \arrow[dd,""] \\
e_{1,2}' \arrow[rr,"" near end] \arrow[dr,swap,""] && e_1
\arrow[dr,swap,""]
\\
& e_2 \arrow[rr,""] \arrow[uu,<-,crossing over,"" near end]&& m.
\end{tikzcd} 
\end{equation*}
We will need the following relations for the core elements of the lower faces
in detailed form. 
\begin{gather}
\label{eq:w23-detailed}
e_{2,3}\subtraction\limits_{E_2}e_{2,3}' = 
\tilde{0}^{2,3}_{e_2}\addition\limits_{E_3} w_{23},\quad
e_{2,3}\subtraction\limits_{E_3}e_{2,3}' = 
\tilde{0}^{2,3}_{e_3}\addition\limits_{E_2} w_{23},\\
\label{eq:e13}
e_{1,3}'\subtraction\limits_{E_1}e_{1,3} =
\tilde{0}^{1,3}_{e_1}\addition\limits_{E_3}w_{13},\quad
e_{1,3}'\subtraction\limits_{E_3}e_{1,3} =
\tilde{0}^{1,3}_{e_3}\addition\limits_{E_1}w_{13},\\
\label{eq:e12}
e_{1,2}\subtraction\limits_{E_1}e_{1,2}' =
\tilde{0}^{1,2}_{e_1}\addition\limits_{E_2} w_{12},
\quad
e_{1,2}\subtraction\limits_{E_2}e_{1,2}' =
\tilde{0}^{1,2}_{e_2}\addition\limits_{E_1} w_{12},
\end{gather}
where $w_{23}\in E_{23}$, $w_{13}\in E_{13}$ and 
$w_{12}\in E_{12}$.

For the zeros of these $w$ elements, the diagrams are
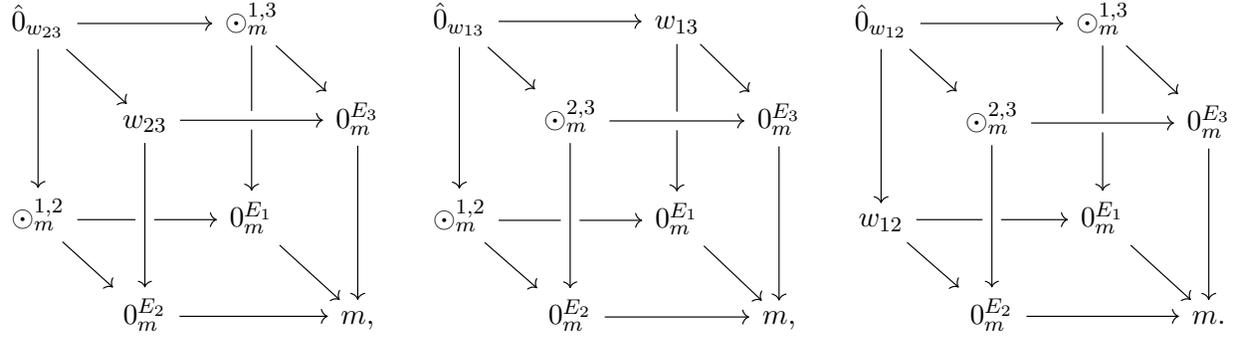
\begin{equation*}
\begin{tikzcd}[row sep=1.5em, column sep=1.15em]
\hat{0}_{w_{23}} \arrow[rr] \arrow[dr,swap,""]
\arrow[dd,swap,""] && \odot^{1,3}_m \arrow[dd,swap,"" near start] \arrow[dr,""]
\\
& w_{23}\arrow[rr,crossing over,"" near start] &&
  0^{E_3}_m \arrow[dd,""] \\
  \odot^{1,2}_m \arrow[rr,"" near end] \arrow[dr,swap,""] &&
  0^{E_1}_m
\arrow[dr,swap,""]
\\
& 0^{E_2}_m \arrow[rr,""] \arrow[uu,<-,crossing over,"" near end]&& m,
\end{tikzcd}
\quad
\begin{tikzcd}[row sep=1.5em, column sep=1.15em]
\hat{0}_{w_{13}} \arrow[rr] \arrow[dr,swap,""]
\arrow[dd,swap,""] && w_{13} \arrow[dd,swap,"" near start]
\arrow[dr,""]
\\
&\odot^{2,3}_m \arrow[rr,crossing over,"" near start] &&
  0^{E_3}_m \arrow[dd,""] \\
  \odot^{1,2}_m \arrow[rr,"" near end] \arrow[dr,swap,""] && 0^{E_1}_m
\arrow[dr,swap,""]
\\
& 0^{E_2}_m \arrow[rr,""] \arrow[uu,<-,crossing over,"" near end]&& m,
\end{tikzcd}
\quad
\begin{tikzcd}[row sep=1.5em, column sep=1.15em]
\hat{0}_{w_{12}} \arrow[rr] \arrow[dr,swap,""]
\arrow[dd,swap,""] && \odot^{1,3}_m \arrow[dd,swap,"" near start] \arrow[dr,""]
\\
&\odot^{2,3}_m \arrow[rr,crossing over,"" near start] &&
  0^{E_3}_m \arrow[dd,""] \\
 w_{12}\arrow[rr,"" near end] \arrow[dr,swap,""] && 0^{E_1}_m
\arrow[dr,swap,""]
\\
& 0^{E_2}_m \arrow[rr,""] \arrow[uu,<-,crossing over,"" near end]&& m.
\end{tikzcd}
\end{equation*}

\subsubsection{Core and ultracore elements arising from the grid}
\label{ssect:differences}

We collect here for reference the definitions and outlines of the
core and ultracore elements arising from the grid. 

\subsubsection{$\bullet$\ $\la_1, k_1$ and $u_1$}

The elements $\one$ and $\two$ have the same Right and Back
faces, and so their differences define an element $\lambda_1\in E_{23,1}$
with outline
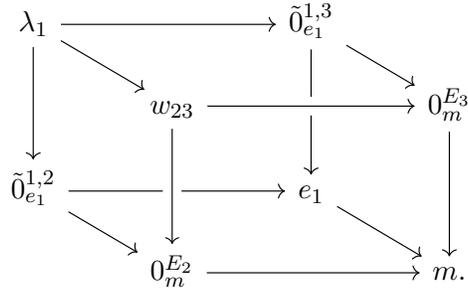
\begin{equation*}
\begin{tikzcd}[row sep=1em]
\la_1 \arrow[rr] \arrow[dr,swap,""] \arrow[dd,swap,""] &&
  \tilde{0}^{1,3}_{e_1} \arrow[dd,swap,"" near start] \arrow[dr,""] \\
& w_{23} \arrow[rr,crossing over,"" near start] &&
 0^{E_3}_m \arrow[dd,""] \\
\tilde{0}^{1,2}_{e_1} \arrow[rr,"" near end] \arrow[dr,swap,""] &&
e_1
\arrow[dr,swap,""]
\\
&0^{E_2}_m\arrow[rr,""] \arrow[uu,<-,crossing over,"" near end]&& m.
\end{tikzcd} 
\end{equation*}
Using (\ref{same_right_and_bottom_faces_differences}) the defining equations are
\begin{equation}
\label{eq:12}
\one\subtraction\limits_{1,2}\two =
\hat{0}_{e_{1,2}}\addition\limits_{1,3}\la_1,\qquad
\one\subtraction\limits_{1,3}\two =
\hat{0}_{e_{1,3}}\addition\limits_{1,2}\la_1.
\end{equation}
If we look at $\thr$ and $\six$, we see that they also have two faces
in common, and their differences define a $k_1\in E_{23,1}$, with
outline
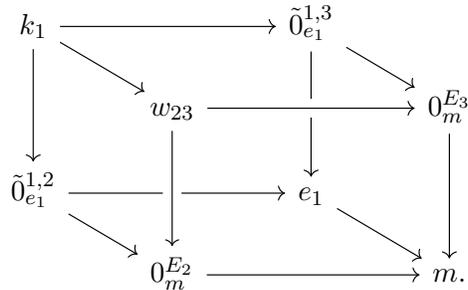
\begin{equation*}
\begin{tikzcd}[row sep=1em]
k_1 \arrow[rr] \arrow[dr,swap,""] \arrow[dd,swap,""] &&
  \tilde{0}^{1,3}_{e_1} \arrow[dd,swap,"" near start] \arrow[dr,""] \\
& w_{23} \arrow[rr,crossing over,"" near start] &&
 0^{E_3}_m \arrow[dd,""] \\
\tilde{0}^{1,2}_{e_1} \arrow[rr,"" near end] \arrow[dr,swap,""] &&
e_1
\arrow[dr,swap,""]
\\
&0^{E_2}_m\arrow[rr,""] \arrow[uu,<-,crossing over,"" near end]&& m.
\end{tikzcd} 
\end{equation*}
The two differences are, again using
(\ref{same_right_and_bottom_faces_differences}),
\begin{equation*}
\thr\subtraction\limits_{1,2}\six =
\hat{0}_{e_{1,2}'}\addition\limits_{1,3}k_1,\qquad
\thr\subtraction\limits_{1,3}\six =
\hat{0}_{e_{1,3}'}\addition\limits_{1,2}k_1.
\end{equation*}
We see that $\la_1$ and $k_1$ have the same outlines so they 
differ by an ultracore element $u_1\in E_{123}$,
\begin{subequations}
\begin{eqnarray}
\lambda_1\subtraction\limits_{1,3}k_1 & = &
\hat{0}_{e_1}\addition\limits_{2,3}u_1,\label{eq:u1a}\\
\la_1\subtraction\limits_{1,2}k_1 & = &
\hat{0}_{e_1}\addition\limits_{2,3}u_1,\label{eq:u1b}\\
\lambda_1\subtraction\limits_{2,3}k_1& = &
\hat{0}_{w_{23}}\addition\limits_{1,3}u_1 = 
\hat{0}_{w_{23}}\addition\limits_{1,2}u_1.\label{eq:u1c}
\end{eqnarray}
\end{subequations}
There are four ways of describing the warp
\begin{subequations}
\begin{align}
(\one\subtraction\limits_{1,2}\two)\subtraction\limits_{1,3}(\thr\subtraction\limits_{1,2}\six)
&=\hat{0}_{e_1}\addition\limits_{2,3}
  (\hat{0}_{w_{12}}\addition\limits_{1,3/2,3}u_1), \label{1a}\\
(\one\subtraction\limits_{1,2}\two)\subtraction\limits_{2,3}(\thr\subtraction\limits_{1,2}\six)
&=(\hat{0}_{w_{12}}\addition\limits_{1,3}\hat{0}_{w_{23}})\addition\limits_{1,3}
  (\hat{0}_{e_2}\addition\limits_{1,3}u_1), \label{1b}\\
(\one\subtraction\limits_{1,3}\two)\subtraction\limits_{1,2}(\thr\subtraction\limits_{1,3}\six)
&=\hat{0}_{e_1}\addition\limits_{2,3}
  (\hat{0}_{-w_{13}}\addition\limits_{1,2/2,3}u_1), \label{1c}\\
(\one\subtraction\limits_{1,3}\two)\subtraction\limits_{2,3}(\thr\subtraction\limits_{1,3}\six)
&=(\hat{0}_{-w_{13}}\addition\limits_{1,2}\hat{0}_{w_{23}})
  \addition\limits_{1,2}(\hat{0}_{e_3}\addition\limits_{1,2}u_1). \label{1d}
\end{align}
\end{subequations}

Here we write $\hat{0}_{w_{12}}\addition\limits_{1,3/2,3}u_1$ to denote
$\hat{0}_{w_{12}}\addition\limits_{1,3}u_1=\hat{0}_{w_{12}}\addition\limits_{2,3}u_1$.

\subsubsection{$\bullet$\ $\la_2, k_2$ and $u_2$}

The same procedure can be applied to $\thr$ and $\fou$;
they have the same Front and Down faces, so their differences will define an
element $\la_2\in E_{13,2}$
\begin{equation*}
\begin{tikzcd}[row sep=1em]
\la_2 \arrow[rr] \arrow[dr,swap,""] \arrow[dd,swap,""] &&
  w_{13} \arrow[dd,swap,"" near start] \arrow[dr,""] \\
& \tilde{0}^{2,3}_{e_2} \arrow[rr,crossing over,"" near start] &&
 0^{E_3}_m \arrow[dd,""] \\
\tilde{0}^{1,2}_{e_2} \arrow[rr,"" near end] \arrow[dr,swap,""] &&
0^{E_1}_m
\arrow[dr,swap,""]
\\
&e_2\arrow[rr,""] \arrow[uu,<-,crossing over,"" near end]&& m.
\end{tikzcd} 
\end{equation*}
The corresponding equations, using (\ref{same_front_and_bottom_faces_differences}), are
\begin{equation*}
\thr\subtraction\limits_{1,2}\fou =
\hat{0}_{e_{1,2}'}\addition\limits_{2,3}\la_2,\qquad
\thr\subtraction\limits_{2,3}\fou =
\hat{0}_{e_{2,3}}\addition\limits_{1,2}\la_2.
\end{equation*}

If we look at $\fiv$ and $\two$, their differences define a
$k_2\in E_{13,2}$, with outline
\begin{equation*}
\begin{tikzcd}[row sep=1em]
k_2 \arrow[rr] \arrow[dr,swap,""] \arrow[dd,swap,""] &&
  w_{13} \arrow[dd,swap,"" near start] \arrow[dr,""] \\
& \tilde{0}^{2,3}_{e_2} \arrow[rr,crossing over,"" near start] &&
 0^{E_3}_m \arrow[dd,""] \\
\tilde{0}^{1,2}_{e_2} \arrow[rr,"" near end] \arrow[dr,swap,""] &&
0^{E_1}_m
\arrow[dr,swap,""]
\\
&e_2\arrow[rr,""] \arrow[uu,<-,crossing over,"" near end]&& m,
\end{tikzcd} 
\end{equation*}
and the differences defined are, due to
(\ref{same_front_and_bottom_faces_differences}),
\begin{equation*}
\fiv\subtraction\limits_{1,2}\two =
\hat{0}_{e_{1,2}}\addition\limits_{2,3}k_2,\qquad
\fiv\subtraction\limits_{2,3}\two =
\hat{0}_{e_{2,3}'}\addition\limits_{1,2}k_2.
\end{equation*}
Since $\la_2$ and $k_2$ have the same outlines, they differ by an ultracore
element $u_2\in E_{123}$,
\begin{subequations}
\begin{eqnarray}
\la_2\subtraction\limits_{1,3}k_2 & = &
\hat{0}_{w_{13}}\addition\limits_{1,2/2,3}u_2,\label{eq:u2a}\\
\la_2\subtraction\limits_{1,2}k_2 & = &
\hat{0}_{e_2}\addition\limits_{1,3}u_2,\label{eq:u2b}\\
\la_2\subtraction\limits_{2,3}k_2& = &
\hat{0}_{e_2}\addition\limits_{1,3}u_2.\label{eq:u2c}
\end{eqnarray}
\end{subequations}
Again there are four ways of describing the warp

\begin{subequations}
\begin{align}
(\thr\subtraction\limits_{1,2}\fou)\subtraction\limits_{2,3}(\fiv\subtraction\limits_{1,2}\two)
&=\hat{0}_{e_2}\addition\limits_{1,3}
  (\hat{0}_{-w_{12}}\addition\limits_{1,3/2,3}u_2), \label{2a}\\
(\thr\subtraction\limits_{1,2}\fou)\subtraction\limits_{1,3}(\fiv\subtraction\limits_{1,2}\two)
&=(\hat{0}_{w_{13}}\addition\limits_{2,3}\hat{0}_{-w_{12}})\addition\limits_{2,3}
  (\hat{0}_{e_1}\addition\limits_{2,3}u_2), \label{2b}\\
(\thr\subtraction\limits_{2,3}\fou)\subtraction\limits_{1,2}(\fiv\subtraction\limits_{2,3}\two)
&=\hat{0}_{e_2}\addition\limits_{1,3}
  (\hat{0}_{w_{23}}\addition\limits_{1,2/1,3}u_2), \label{2c}\\
(\thr\subtraction\limits_{2,3}\fou)\subtraction\limits_{1,3}(\fiv\subtraction\limits_{2,3}\two)
&=(\hat{0}_{w_{23}}\addition\limits_{1,2}\hat{0}_{w_{13}})
  \addition\limits_{1,2}(\hat{0}_{e_3}\addition\limits_{1,2}u_2). \label{2d}
\end{align}
\end{subequations}

\subsubsection{$\bullet$\ $\la_3,k_3$ and $u_3$}

Likewise $\fiv$ and $\six$ define $\la_3\in E_{12,3}$ with outline
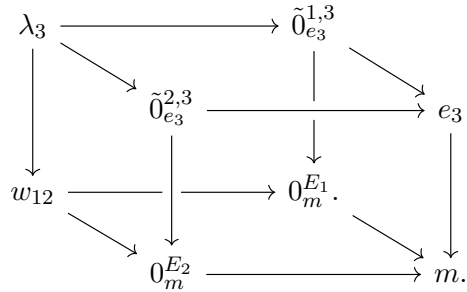
\begin{equation*}
\begin{tikzcd}[row sep=1em]
\la_3 \arrow[rr] \arrow[dr,swap,""] \arrow[dd,swap,""] &&
  \tilde{0}^{1,3}_{e_3} \arrow[dd,swap,"" near start] \arrow[dr,""] \\
& \tilde{0}^{2,3}_{e_3} \arrow[rr,crossing over,"" near start] &&
  e_3 \arrow[dd,""] \\
w_{12} \arrow[rr,"" near end] \arrow[dr,swap,""] && 0^{E_1}_m.
\arrow[dr,swap,""]
\\
& 0^{E_2}_m \arrow[rr,""] \arrow[uu,<-,crossing over,"" near end]&& m.
\end{tikzcd} 
\end{equation*}
The corresponding relations are, due to
(\ref{same_front_and_right_faces_differences}),
\begin{equation*}
\fiv\subtraction\limits_{1,3}\six =
\hat{0}_{e_{1,3}'}\addition\limits_{2,3}\la_3,\qquad
\fiv\subtraction\limits_{2,3}\six =
\hat{0}_{e_{2,3}'}\addition\limits_{1,3}\la_3.
\end{equation*}
Likewise $\one$ and $\fou$ define a $k_3\in E_{12,3}$ with outline
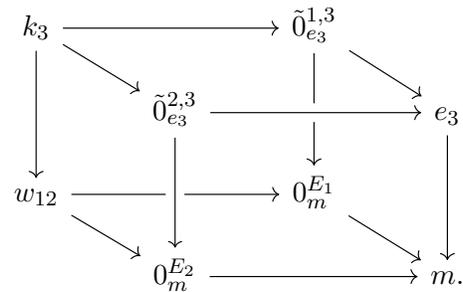
\begin{equation*}
\begin{tikzcd}[row sep=1em]
k_3 \arrow[rr] \arrow[dr,swap,""] \arrow[dd,swap,""] &&
  \tilde{0}^{1,3}_{e_3} \arrow[dd,swap,"" near start] \arrow[dr,""] \\
& \tilde{0}^{2,3}_{e_3} \arrow[rr,crossing over,"" near start] &&
  e_3 \arrow[dd,""] \\
w_{12} \arrow[rr,"" near end] \arrow[dr,swap,""] && 0^{E_1}_m
\arrow[dr,swap,""]
\\
& 0^{E_2}_m \arrow[rr,""] \arrow[uu,<-,crossing over,"" near end]&& m.
\end{tikzcd} 
\end{equation*}
The differences defined are, due to
(\ref{same_front_and_right_faces_differences}),
\begin{equation*}
\one\subtraction\limits_{1,3}\fou =
\hat{0}_{e_{1,3}}\addition\limits_{2,3}k_3,\qquad
\one\subtraction\limits_{2,3}\fou =
\hat{0}_{e_{2,3}}\addition\limits_{1,3}k_3.
\end{equation*}
The ultracore element $u_3\in E_{123}$ defined by $\la_3$ and $k_3$ satisfies
\begin{subequations}
\begin{eqnarray}
\la_3\subtraction\limits_{1,3}k_3 & = & 
\hat{0}_{e_3}\addition\limits_{1,2}u_3,\label{eq:u3a}\\
\la_3\subtraction\limits_{1,2}k_3 &= & 
\hat{0}_{w_{12}}\addition\limits_{1,3/2,3}u_3,\label{eq:u3b}
\\
\la_3\subtraction\limits_{2,3}k_3 & = & 
\hat{0}_{e_3}\addition\limits_{1,2}u_3.\label{eq:u3c}
\end{eqnarray}
\end{subequations}
The four relations in this case are

\begin{subequations}
\begin{align}
(\fiv\subtraction\limits_{1,3}\six)\subtraction\limits_{2,3}(\one\subtraction\limits_{1,3}\fou)
&=\hat{0}_{e_3}\addition\limits_{1,2}
  (\hat{0}_{w_{13}}\addition\limits_{1,2/2,3}u_3), \label{3a}\\
(\fiv\subtraction\limits_{1,3}\six)\subtraction\limits_{1,2}(\one\subtraction\limits_{1,3}\fou)
&=(\hat{0}_{w_{13}}\addition\limits_{2,3}\hat{0}_{w_{12}})\addition\limits_{2,3}
  (\hat{0}_{e_1}\addition\limits_{2,3}u_3), \label{3b}\\
(\fiv\subtraction\limits_{2,3}\six)\subtraction\limits_{1,3}(\one\subtraction\limits_{2,3}\fou)
&=\hat{0}_{e_3}\addition\limits_{1,2}
  (\hat{0}_{-w_{23}}\addition\limits_{1,2/1,3}u_3), \label{3c}\\
(\fiv\subtraction\limits_{2,3}\six)\subtraction\limits_{1,2}(\one\subtraction\limits_{2,3}\fou)
&=(\hat{0}_{-w_{23}}\addition\limits_{1,3}\hat{0}_{w_{12}})
  \addition\limits_{1,3}(\hat{0}_{e_2}\addition\limits_{1,3}u_3). \label{3d}
\end{align}
\end{subequations}

\subsection{Proof of the warp theorem} 

We will show that $u_1+u_2+u_3 = \odot^3_m$ 
by showing that $u_1 = -u_2-u_3$. There
are four steps. 

\smallskip

\textbf{Step 1.} Rewrite (\ref{1b})
\begin{equation*}
(\one\subtraction\limits_{1,2
}\two)\subtraction\limits_{2,3}(\thr\subtraction\limits_{1,2}\six)
\end{equation*}
as
\begin{equation*}
(\one\subtraction\limits_{2,3}
\thr)\subtraction\limits_{1,2}(\two\subtraction\limits_{2,3}\six),
\end{equation*}
using the interchange law for the \dvb formed by the Left face.
We know from (\ref{1b}) that the ultracore element defined by the first expression is $u_1$, 
therefore, the ultracore element of the latter expression will also be $u_1$. 
We will show that the second expression 
has $-u_2-u_3$ as its ultracore element, and this will show that $u_1 = -u_2-u_3$.

\smallskip

\textbf{Step 2.} First, write $\one\subtraction\limits_{2,3}\thr$ as
\begin{equation*}
\one\subtraction\limits_{2,3}\thr = 
(\one\subtraction\limits_{2,3}\fou)\subtraction\limits_{2,3}
(\thr\subtraction\limits_{2,3}\fou),
\end{equation*}
where we have that 
\begin{equation}
\label{eq:1-4,3-4}
\one\subtraction\limits_{2,3}\fou =
\hat{0}_{e_{2,3}}\addition\limits_{1,3}k_3,\qquad
\thr\subtraction\limits_{2,3}\fou =
\hat{0}_{e_{2,3}}\addition\limits_{1,2}\la_2.
\end{equation}

\smallskip

\textbf{Step 3.} Similarly, write $\two\subtraction\limits_{2,3}\six$ as
\begin{equation}
\label{eq:5-6,5-2}
\two\subtraction\limits_{2,3}\six = 
(\fiv\subtraction\limits_{2,3}\six)\subtraction\limits_{2,3}
(\fiv\subtraction\limits_{2,3}\two),
\end{equation}
and we have
\begin{equation*}
\fiv\subtraction\limits_{2,3}\six =
\hat{0}_{e_{2,3}'}\addition\limits_{1,3}\la_3,\qquad
\fiv\subtraction\limits_{2,3}\two =
\hat{0}_{e_{2,3}'}\addition\limits_{1,2}k_2.
\end{equation*}

\smallskip

\textbf{Step 4.} 
Since our convention is that $\la_3-k_3$ defines $u_3$, it follows
that $k_3-\la_3$ defines $-u_3$. 
It is essential to follow these conventions.

We are finally able to complete the proof of Theorem \ref{thm:grid}. 
\begin{eqnarray*}
(\one\subtraction\limits_{2,3}
\thr)\subtraction\limits_{1,2}(\two\subtraction\limits_{2,3}\six)&=&
[(\one\subtraction\limits_{2,3}\fou)\subtraction\limits_{2,3}
(\thr\subtraction\limits_{2,3}\fou)]\subtraction\limits_{1,2}[(\fiv\subtraction\limits_{2,3}\six)\subtraction\limits_{2,3}
(\fiv\subtraction\limits_{2,3}\two)]\\ 
&\phantom{.}& \text{(using operations in single \vbs)}\\
&=&
[(\hat{0}_{e_{2,3}}\addition\limits_{1,3}k_3)\subtraction\limits_{2,3}
(\hat{0}_{e_{2,3}}\addition\limits_{1,2}\la_2)]\subtraction\limits_{1,2}
[(\hat{0}_{e_{2,3}'}\addition\limits_{1,3}\la_3)\subtraction\limits_{2,3}
(\hat{0}_{e_{2,3}'}\addition\limits_{1,2}k_2)]\\
&\phantom{.}& \text{(using (\ref{eq:1-4,3-4}) and (\ref{eq:5-6,5-2}))}\\
& = &
[(\hat{0}_{e_{2,3}}\addition\limits_{1,3}k_3)\subtraction\limits_{1,2}
(\hat{0}_{e_{2,3}'}\addition\limits_{1,3}\la_3)]\subtraction\limits_{2,3}
[(\hat{0}_{e_{2,3}}\addition\limits_{1,2}\la_2)\subtraction\limits_{1,2}
(\hat{0}_{e_{2,3}'}\addition\limits_{1,2}k_2)]\\
&\phantom{.}& \text{(interchange law in Left face applied to outer operations)}\\
& = &
[(\hat{0}_{e_{2,3}}\subtraction\limits_{1,2}\hat{0}_{e_{2,3}'})\addition\limits_{1,3}
(k_3\subtraction\limits_{1,2}\la_3)]\subtraction\limits_{2,3}
[(\hat{0}_{e_{2,3}}\subtraction\limits_{1,2}\hat{0}_{e_{2,3}'})\addition\limits_{1,2}
(\la_2\subtraction\limits_{1,2}k_2)]\\
&\phantom{.}& \text{(interchange law in Back face for first term)}\\
& = &
[(\hat{0}_{w_{23}}\addition\limits_{1,3}\hat{0}_{e_2})\addition\limits_{1,3}
(\hat{0}_{w_{12}}\subtraction\limits_{1,3}u_3)]\subtraction\limits_{2,3}
[(\hat{0}_{w_{23}}\addition\limits_{1,3}\hat{0}_{e_2})\addition\limits_{1,2}
(\hat{0}_{e_2}\addition\limits_{1,3}u_2)]\\
&\phantom{.}& \text{(first term in each $[~]$ by (\ref{eq:w23-2}), then use (\ref{eq:u3b}) and 
(\ref{eq:u2b}))}\\
& = &
[\hat{0}_{w_{23}}\addition\limits_{1,3}(\hat{0}_{e_2}\addition\limits_{1,3}
\hat{0}_{w_{12}}\subtraction\limits_{1,3}u_3)]\subtraction\limits_{2,3}
[(\hat{0}_{w_{23}}\addition\limits_{1,2}u_2)
\addition\limits_{1,3}(\hat{0}_{e_2}\addition\limits_{1,2}\hat{0}_{e_2})]\\
&\phantom{.}& \text{(second group by interchange law in Back face)}\\
& = &
[\hat{0}_{w_{23}}\addition\limits_{1,3}(\hat{0}_{e_2}\addition\limits_{1,3}
\hat{0}_{w_{12}}\subtraction\limits_{1,3}u_3)]\subtraction\limits_{2,3}
[(\hat{0}_{w_{23}}\addition\limits_{1,3}u_2)
\addition\limits_{1,3}\hat{0}_{e_2}]\\
&\phantom{.}& \text{(zeros in final group are over base of the addition)}\\
& = &
[\hat{0}_{w_{23}}\addition\limits_{1,3}(\hat{0}_{e_2}\addition\limits_{1,3}
\hat{0}_{w_{12}}\subtraction\limits_{1,3}u_3)]\subtraction\limits_{2,3}
[\hat{0}_{w_{23}}\addition\limits_{1,3}(\hat{0}_{e_2}\addition\limits_{1,3}u_2)]\\
&\phantom{.}& \text{(second group is in ordinary \vb)}\\
& = &
[\hat{0}_{w_{23}}\subtraction\limits_{2,3}\hat{0}_{w_{23}}]\addition\limits_{1,3}
[(\hat{0}_{e_2}\addition\limits_{1,3}\hat{0}_{w_{12}}\subtraction\limits_{1,3}
u_3)\subtraction\limits_{2,3}(\hat{0}_{e_2}\addition\limits_{1,3}u_2)]\\
&\phantom{.}& \text{(interchange law in Up face)}\\
& = &
\hat{0}_{w_{23}}\addition\limits_{1,3}[\hat{0}_{e_2}\subtraction\limits_{2,3}\hat{0}_{e_2}]
\addition\limits_{1,3}[(\hat{0}_{w_{12}}\subtraction\limits_{1,3}
u_3)\subtraction\limits_{2,3}u_2]\\
&\phantom{.}& \text{(zeros are zeros over base of addition; then interchange law in Up face)}\\
& = &
\hat{0}_{w_{23}}\addition\limits_{1,3}\hat{0}_{e_2}\addition\limits_{1,3}
[(\hat{0}_{w_{12}}\subtraction\limits_{2,3}
u_3)\subtraction\limits_{2,3}u_2]\\
&\phantom{.}& \text{(zeros are zeros over base of addition)}\\
& = &
\hat{0}_{w_{23}}\addition\limits_{1,3}\hat{0}_{e_2}\addition\limits_{1,3}
[\hat{0}_{w_{12}}\subtraction\limits_{2,3}
(u_3\addition\limits_{2,3}u_2)]\\
& = &
\hat{0}_{w_{23}}\addition\limits_{1,3}\hat{0}_{e_2}\addition\limits_{1,3}
\hat{0}_{w_{12}}\subtraction\limits_{1,3}(u_3\addition\limits_{1,3}u_2)\\
&\phantom{.}& \text{(using an equation of the form (\ref{eq:13,12}))}\\
\end{eqnarray*}
from which we obtain $-(u_3+u_2)$ as the ultracore element.

Comparing this with (\ref{1b}),
$$
(\one\subtraction\limits_{1,2}\two)\subtraction\limits_{2,3}(\thr\subtraction\limits_{1,2}\six)
= 
(\hat{0}_{w_{12}}\addition\limits_{1,3}\hat{0}_{w_{23}})
\addition\limits_{1,3}  
(\hat{0}_{e_2}\addition\limits_{1,3}u_1),
$$
we have $u_1 = -(u_3+u_2)$ as desired. 

This completes the proof of the warp theorem. 

The strategy of this proof deserves some commentary. 

What should the warp of a grid on a \tvb be? Or, in other words, why are we
interested in the ultrawarps of a grid of a \tvb?

The warp of a grid in the double case is a section of the core \vb, and
measures the non-commutativity of the two routes defined by the grid.

So far, we have seen that all operations on a \tvb are iterations of operations
defined in \dvbs. The ultracore, for example, is the core of the core \dvbs.

For these reasons, we would want the warp of a grid 
in the triple case to be a section of the ultracore \vb, and to measure the
non-commutativity of routes defined by the grid.

Pick an upper face of $E$, for example the Up face. If we compare the two
routes defined by the grid in this face, then we obtain an element of the
(U-D) core \dvb, which we denoted by $\lambda_3$. Similarly for the other
upper faces, the non-commutativity of the corresponding routes defines $\la_1$
and $\la_2$. The three $\lambda$'s are elements of different spaces; therefore,
if we tried to compare them, or indeed perform any sort of operation with them
(such as adding them or subtracting them), we would see that such an operation
could be algebraically possible but would not be geometrically meaningful.

The same applies for the three $k_i$ defined by the comparison of the routes for
the lower faces.

The $\la_i$'s and the corresponding $k_i$'s however, are elements of the same
spaces, therefore, comparing them is a possibility, and indeed the only
sensible operation. And by comparing them, we measure the
non-commutativity of four routes, instead of two.

This can be done for the three pairs of $\la_i$ and $k_i$, and so 
we obtain the three ultrawarps.

So what does the warp theorem tell us?

Each ultrawarp measures the non-commutativity of four routes. In total, a grid on
a \tvb provides six different routes from $M$ to $E$. The sum of the three
ultrawarps takes into account each route twice, once with a positive and once
with a negative sign. The warp theorem tells us that these add up to zero,
a result that seems reasonable. The different \vb structures over which the
operations take place however, are the main obstacle here --- as soon as
one realizes that simple operations like addition and subtraction in the \tvb
setting are no longer simple.

\newpage

\section{The \tvb $T^2A$ and connections in $A$}  
\label{sect:T2E}

In this section and the next we examine two typical 
instances of the warp theorem. In this section we consider the \tvb $T^2A$
where $A\to M$ is a vector bundle, and grids which arise in it from
connections in $A$. In the following section we consider $T^3M$, the
triple tangent bundle of a manifold $M$. 

For a vector bundle $(A, q, M)$ there is a \tvb structure on $T^2A$ as shown in
(\ref{eq:T2A}). 
\begin{equation}
  \label{eq:T2A}
\begin{tikzcd}[row sep=2.5em, column sep=2.5em]
T^2A \arrow[rr,"T^2(q)"] \arrow[dr,"T(p_A)"] \arrow[dd,swap,"p_{TA}"] &&
 T^2M \arrow[dd,swap,"p_{TM}" near start] \arrow[dr,"T(p)"] \\
& TA \arrow[rr,crossing over,"T(q)" near start] &&
  TM \arrow[dd,"p"] \\
TA \arrow[rr,"T(q)" near end] \arrow[dr,swap,"p_A"] && TM
\arrow[dr,swap,"p"]
\\
& A \arrow[rr,"q"] \arrow[uu,<-,crossing over,"p_A" near end]&& M.
\end{tikzcd}
\end{equation}
Here the Down face is the usual \dvb $TA$ and the Up face is the tangent
prolongation of this; that is, it is obtained by applying the tangent functor
to each structure in $TA$. Each vertical vector bundle in (\ref{eq:T2A})
is a standard tangent bundle. 

\subsection{The core \dvbs of $T^2A$}

The three core \dvbs are shown in (\ref{eq:T2Acores}), in the usual order
(B-F), (L-R), and (U-D), 
and arranged 
as in (\ref{eq:UP-DOWN CORE DVB}). 
\begin{equation}
\label{eq:T2Acores}
\begin{tikzcd}[row sep=1.5cm, column sep = 1.5cm]
TA \arrow[r,"p_A"]\arrow[d, swap,"T(q)"] 
&A \arrow[d,"q"]\\
TM\arrow[r,swap,"p"]
&M, 
\end{tikzcd}
\qquad
\begin{tikzcd}[row sep=1.5cm, column sep = 1.5cm]
TA \arrow[r,"T(q)"]\arrow[d, swap,"p_A"] 
&TM \arrow[d,"p"]\\
A\arrow[r,swap,"q"]
&M, 
\end{tikzcd} 
\qquad
\begin{tikzcd}[row sep=1.5cm, column sep = 1.5cm]
TA \arrow[r,"p_A"]\arrow[d, swap,"T(q)"] 
&A \arrow[d,"q"]\\
TM\arrow[r,swap,"p"]
&M. 
\end{tikzcd}
\end{equation}

These core \dvbs are the same as abstract \dvbs but are embedded differently in $T^2A$. 

An element $\xi\in TA$ determines core elements of the Back, Left and Up faces. 
We denote these by $\Bar{\xi}^{B}$, $\Bar{\xi}^{L}$ and $\Bar{\xi}^{U}$,
respectively.
Since they are elements of $T(TA)$ they can be represented as tangent vectors to curves
in $TA$. 

The Back face is the tangent \dvb for the tangent prolongation 
bundle $T(q)\co TA\to TM$. For the core of the Back face we have 
\begin{equation}
  \label{eq:coreB}
\Bar{\xi}^{B} = \frac{d}{dt}(t\mult\limits_{TM}\xi)\Big\rvert_{t = 0}
\end{equation}
where $t\mult\limits_{TM}\xi$ denotes scalar multiplication in 
$T(q)\co TA\to TM$. 

The Left face is the double tangent bundle of the manifold $A$. Given $\xi\in TA$
the core element is
\begin{equation}
    \label{eq:coreL}
\Bar{\xi}^{L} = \frac{d}{dt}(t\xi)\Big\rvert_{t = 0}
\end{equation}
where the scalar multiplication is in the usual tangent bundle $TA\to A$. 

For the core of the Up face, first write $\xi = \left.\frac{d}{dt}a_t\right|_{t=0}$
where $a_t$ is a curve in $A$. Write $\Bar{a_t}\in TA$ for the core element
corresponding to $a_t$. Then 
\begin{equation}
    \label{eq:coreU}
\Bar{\xi}^{U} = \frac{d}{dt}(\Bar{a_t})\Big\rvert_{t = 0}. 
\end{equation}

\subsection{The canonical involution on $T^2A$}

The canonical involution $J_A\co T^2A\to T^2A$ for the manifold $A$
is an isomorphism from the \dvb $T^2A$ to its flip. In what follows
we will need to use it as a map of \tvbs.  

\begin{proposition}
\label{prop:JAmorphism}
The map $J_A$ is an isomorphism of \tvbs as shown in {\rm (\ref{J_A})}. 
\end{proposition}

\begin{equation}\label{J_A}
\begin{tikzcd}[row sep=2.5em, column sep=2.5em]
T^2A \arrow[rr,"T^2(q)"] \arrow[dr,"T(p_A)",near end] \arrow[dd,swap,"p_{TA}"]
&& T^2M \arrow[dd,swap,"p_{TM}" near start] \arrow[dr,"T(p)"] \\
& TA \arrow[rr,crossing over,"T(q)" near start] &&
  TM \arrow[dd,"p"] \\
TA \arrow[rr,"T(q)" near end] \arrow[dr,swap,"p_A"] && TM
\arrow[dr,swap,"p"]
\\
& A \arrow[rr,"q"] \arrow[uu,<-,crossing over,"p_A" near end]&& M
\end{tikzcd}
\xRightarrow{J_A}
\begin{tikzcd}[row sep=2.5em, column sep=2.5em]
T^2A \arrow[rr,"T^2(q)"] \arrow[dr,"p_{TA}",near end] \arrow[dd,swap,"T(p_A)"]
&& T^2M \arrow[dd,swap,"T(p)" near start] \arrow[dr,"p_{TM}"] \\
& TA \arrow[rr,crossing over,"T(q)" near start] &&
  TM \arrow[dd,"p"] \\
TA \arrow[rr,"T(q)" near end] \arrow[dr,swap,"p_A"] && TM
\arrow[dr,swap,"p"]
\\
& A \arrow[rr,"q"] \arrow[uu,<-,crossing over,"p_A" near end]&& M.
\end{tikzcd}
\end{equation} 

In (\ref{J_A}) the Left faces are the double tangent bundles of the manifold
$A$ and $J_A$ maps the Left face of the domain to its flip. It interchanges
the Up and Back faces. The Right faces are the double tangent bundles of 
$M$ and $J_A$ induces $J_M\co T^2M\to T^2M$ which maps the Right face of
the domain to its flip. The Front and Down faces are interchanged. 

The proof of Proposition \ref{prop:JAmorphism} relies on 
the following two lemmas. The first is the naturality property of the
canonical involution. 

\begin{lemma}
\label{lem:Jnat}
Let $M$ and $N$ be smooth manifolds, and $F: M\rightarrow N$ a smooth
map. Then $T^2(F)\circ J_M = J_N\circ T^2(F)$, where $T^2(F) = T(T(F))$ is the
tangent of the tangent map $T(F)$.
\end{lemma}

\begin{lemma}
\label{lem:Jadd}
Given $\Phi_1,\Phi_2\in T^2A$, over the same $\xi\in T^2M$, we have
\begin{equation*}
J_A\left(\Phi_1\addition\limits_{T^2(q)} \Phi_2 \right) =
J_A(\Phi_1)\addition\limits_{T^2(q)} J_A(\Phi_2).
\end{equation*}
\end{lemma}

Consider now the maps which $J_A$ induces on the cores. 

Take an element $\xi\in TA$ in the core of the Back face. Regarded
as an element of $T^2A$ this is ${\overline{\xi}}^{B}$, with
outline shown on the left of (\ref{triple-outline-of-(B-F)-core-dvb}).

\begin{equation}\label{triple-outline-of-(B-F)-core-dvb}
\begin{tikzcd}[row sep=1.5em, column sep=1.5em]
\overline{\xi}^{B}\arrow[rr,"T^2(q)"] \arrow[dr,swap,"T(p_A)",near end]
\arrow[dd,swap,"p_{TA}"] && 0^{T^2M}_x \arrow[dd,swap," " near start]
\arrow[dr," "]
\\
& \overline{a} \arrow[rr,crossing over," " near start] &&
  0^{TM}_m \arrow[dd," "] \\
T(0^A)(x) \arrow[rr," " near end] \arrow[dr,swap," "] && x
\arrow[dr,swap," "]
\\
& 0^A_m \arrow[rr," "] \arrow[uu,<-,crossing over," " near end]&& m,
\end{tikzcd}
\qquad
\begin{tikzcd}[row sep=1.5em, column sep=1.5em]
\overline{\xi}^{U}\arrow[rr,"T^2(q)"] \arrow[dr,swap,"T(p_A)",near end]
\arrow[dd,swap,"p_{TA}"] && T(0^{TM})(x) \arrow[dd,swap," " near start]
\arrow[dr," "]
\\
& T(0^A)(x) \arrow[rr,crossing over," " near start] &&
  x \arrow[dd," "] \\
\overline{a} \arrow[rr," " near end] \arrow[dr,swap," "] && 0^{TM}_m
\arrow[dr,swap," "]
\\
& 0^A_m \arrow[rr," "] \arrow[uu,<-,crossing over," " near end]&& m.
\end{tikzcd}
\end{equation}

It follows from (\ref{eq:coreB}) and (\ref{eq:coreU}) that
\begin{equation}\label{J_A:(B-F)to(U-D)}
J_A({\overline{\xi}}^{B}) = {\overline{\xi}}^{U}.
\end{equation}

Since $J_A^2$ is the identity, we also have
\begin{equation}\label{J_A:(U-D)to(B-F)}
J_A({\overline{\xi}}^{U}) = {\overline{\xi}}^{B}. 
\end{equation}

Since the Left faces in (\ref{J_A}) are the double tangent bundle $T^2A$,
the map on the cores of the Left faces is the identity and so
\begin{equation}\label{J_A:LtoL}
J_A({\overline{\xi}}^{L}) = {\overline{\xi}}^{L}. 
\end{equation}

\subsection{Grids on $T^2A$}

Now consider a connection $\nabla$ in $A$. Recall that Example~\ref{ex:conn} 
gave a construction of a grid in $TA$ for which the warp is $\nabla_X\mu$. 
We now extend this idea to define a grid in $T^2A$.

Let $X,Z\in \X(M)$, and $\mu\in \Gamma A$. 
Define the following three double linear sections:
\begin{itemize}
  \renewcommand{\itemsep}{-2pt} 
\item From Front to Back face: $(T(X^{\hlift}); X^{\hlift}, T(X);X)$.
\item From Right to Left face: $(T^2(\mu);T(\mu),T(\mu);\mu)$.
\item From Down to Up face: $(\wtilde{Z^{\hlift}}^A;
Z^{\hlift},\widetilde{Z};Z)$.
\end{itemize}
Here $\wtilde{Z} = J_M\circ T(Z)$ is the complete (or tangent) lift of $Z$ to a
vector field on $TM$. Likewise $\wtilde{Z^{\hlift}}^A$ is the complete lift of
$Z^{\hlift}\in\X(A)$ to a vector field on $TA$. The grid is shown in
(\ref{eq:T2Agrid}).

\begin{equation}
  \label{eq:T2Agrid}
\begin{tikzcd}[row sep=2.5em, column sep=2.5em]
T^2A \arrow[rr,<-,"T^2(\mu)"] \arrow[dr,<-,"T(X^{\hlift})"]
\arrow[dd,<-,swap,"\wtilde{Z^{\hlift}}^A"] && T^2M \arrow[dd,<-,swap,"\widetilde{Z}" near start]
\arrow[dr,<-,"T(X)"] \\
& TA \arrow[rr,<-,crossing over,"T(\mu)" near start] &&
  TM \arrow[dd,<-,"Z"] \\
TA \arrow[rr,<-,"T(\mu)" near start] \arrow[dr,<-,swap,"X^{\hlift}"] && TM
\arrow[dr,<-,"X"]
\\
& A \arrow[rr,<-,swap,"\mu"] \arrow[uu,crossing over,swap,"Z^{\hlift}" near end]&& M. 
\end{tikzcd}
\end{equation}

The core morphisms of the linear double sections will be needed later:
\begin{itemize}
\renewcommand{\itemsep}{-2pt} 
\item For $(T(X^{\hlift}); X^{\hlift}, T(X);X)$ the core morphism is
$(X^{\hlift}, X)$.
\item For $(T^2(\mu);T(\mu),T(\mu);\mu)$ the core morphism is $(T(\mu),
\mu)$.
\item For $(\wtilde{Z^{\hlift}}^A;Z^{\hlift},\widetilde{Z};Z)$
 the core morphism is $(Z^{\hlift},Z)$.
\end{itemize}

The first two cases are instances of the general fact that given a morphism $(\ph,f)$
of vector bundles, the core morphism of the \dvb map $(T(\ph);T(f),\ph;f)$ is $(\ph,f)$. 

To calculate the core morphism of
$(\wtilde{Z^{\hlift}}^A;Z^{\hlift},\widetilde{Z};Z)$, focus on (\ref{eq:Zlds}).
At this point we investigate this linear double section further; it is a \dvb
morphism from the Down face to the Up face of $T^2A$. Note that (\ref{eq:Zlds})
is not a \tvb.
\begin{equation}
  \label{eq:Zlds}
\begin{tikzcd}[row sep=2.5em, column sep=2.5em]
T^2A \arrow[rr,"T^2(q)"] \arrow[dr,"T(p_A)"]
\arrow[dd,<-,swap,"\wtilde{Z^{\hlift}}^A",red] && T^2M
\arrow[dd,<-,swap,"\widetilde{Z}" near start,red] \arrow[dr,"T(p)"] \\
& TA \arrow[rr,crossing over,"T(q)" near start] &&
  TM \arrow[dd,<-,"Z",red] \\
TA \arrow[rr,"T(q)" near start] \arrow[dr,swap,"p_A"] && TM
\arrow[dr,"p"]
\\
& A \arrow[rr,swap,"q"] \arrow[uu,crossing over,swap,"Z^{\hlift}" near
end,red]&& M.
\end{tikzcd}
\end{equation}

Take an element $a\in A$. As an element of the core of the Down face of (\ref{eq:T2Agrid})
it is $\overline{a} = \frac{d}{dt}ta\Big\rvert_{t = 0}\in TA$.
Using the fact that $(\wtilde{Z},Z)$ is a \vb map, we have that
$\wtilde{Z}(0^{TM}_m)= T(0^{TM})(Z(m))$. 

Similarly,
using the fact that $(Z^{\hlift},Z)$ is a \vb map, we have that
$Z^{\hlift}(0^A_m)= T(0^A)(Z(m))$. 

Finally,
\begin{equation*}
\wtilde{Z^{\hlift}}^A(\overline{a}) =
J_A\left(T(Z^{\hlift})(\overline{a})\right) = 
J_A\left(\frac{d}{dt}Z^{\hlift}(ta)\Big\rvert_{t = 0} \right)=
J_A\left(\frac{d}{dt}tZ^{\hlift}(a)\Big\rvert_{t = 0} \right) 
= J_A\left(\overline{Z^{\hlift}(a)}^{B} \right).
\end{equation*}

Note the following. Initially, $\overline{a}\in TA$ is in the core of the Down
face of $T^2A$. The canonical involution $J_A$ maps the Down face to the Front
face (see (\ref{J_A})). Therefore, in $T(Z^{\hlift})(\overline{a})$, $\overline{a}$ is now an
element of the core of the Front face. The maps
$(T(Z^{\hlift});Z^{\hlift},T(Z);Z)$ form a \dvb morphism from the Front to the
Back face of (\ref{eq:T2A}), with core morphism $(Z^{\hlift},Z)$ as usual.
Therefore, $T(Z^{\hlift})(\overline{a}) = \overline{Z^{\hlift}(a)}^{B}$ is now
in the core of the Back face.
And by (\ref{J_A:(B-F)to(U-D)}), it follows that
$J_A\left(\overline{Z^{\hlift}(a)}^{B}\right) = \overline{Z^{\hlift}(a)}^{U}$.
 
This completes the proof that the core morphism of
$(\wtilde{Z^{\hlift}}^A; Z^{\hlift},\widetilde{Z};Z)$ is $(Z^{\hlift},Z)$.

\subsection{The warp of the Back face}
\label{subsect:wbf}

The warp of the Back face is given by
\begin{equation*}
T^2(\mu)(\wtilde{Z}(X(m)))\subtraction\limits_{p_{TA}}
\wtilde{Z^{\hlift}}^A(T(\mu)(X(m))).
\end{equation*}
The outlines of the two elements are
\begin{equation*}
\begin{tikzcd}[row sep=1.5em, column sep=1.5em]
T^2(\mu)(\wtilde{Z}(X(m))) \arrow[rr," "] \arrow[dr,swap," "] \arrow[dd,swap,"
"] && \wtilde{Z}(X(m)) \arrow[dd,swap," " near start] \arrow[dr," "] \\
& T(\mu)(Z(m)) \arrow[rr,crossing over," " near start] &&
  Z(m)\arrow[dd," "] \\
T(\mu)(X(m)) \arrow[rr," " near end] \arrow[dr,swap," "] && X(m)
\arrow[dr,swap," "]
\\
& \mu(m) \arrow[rr," "] \arrow[uu,<-,crossing over," " near end]&& m,
\end{tikzcd}
\end{equation*}
\begin{equation*}
\begin{tikzcd}[row sep=1.5em, column sep=1.5em]
\wtilde{Z^{\hlift}}^A(T(\mu)(X(m))) \arrow[rr," "] \arrow[dr,swap," "]
\arrow[dd,swap," "] && \wtilde{Z}(X(m)) \arrow[dd,swap," " near start] \arrow[dr," "] \\
& Z^{\hlift}(\mu(m)) \arrow[rr,crossing over," " near start] &&
  Z(m)\arrow[dd," "] \\
T(\mu)(X(m)) \arrow[rr," " near end] \arrow[dr,swap," "] && X(m)
\arrow[dr,swap," "]
\\
& \mu(m) \arrow[rr," "] \arrow[uu,<-,crossing over," " near end]&& m,
\end{tikzcd} 
\end{equation*}
(compare with the general triple outlines of the elements $\two$ and $\one$, of
subsection \ref{subsect:notation-of-six-elements}).

Writing the complete lifts as 
$\wtilde{Z^{\hlift}}^A = J_A\circ T(Z^{\hlift})$ and $\wtilde{Z} = J_M\circ
T(Z)$, and using the naturality of $J$-maps (Lemma~\ref{lem:Jnat}), we have that 
\begin{multline}\label{first_line}
T^2(\mu)(\wtilde{Z}(X(m)))\subtraction\limits_{p_{TA}}
\wtilde{Z^{\hlift}}^A(T(\mu)(X(m)))\\
 = T^2(\mu)\left(J_M(T(Z)(X(m)))\right)
\subtraction\limits_{p_{TA}}
J_A\left(T(Z^{\hlift})(T(\mu)(X(m)))\right)\\
 = J_A(T^2(\mu)(T(Z)(X(m))))
\subtraction\limits_{p_{TA}}
J_A(T(Z^{\hlift})(T(\mu)(X(m))))\\
 = J_A (T^2(\mu)(T(Z)(X(m)))) 
\subtraction\limits_{p_{TA}}
J_A (T(Z^{\hlift})(T(\mu)(X(m)))).
\end{multline}

Since $J_A$ interchanges the structures $p_{TA}$ and $T(p_A)$, we can
rewrite the last expression in (\ref{first_line}) as
\begin{equation*}
J_A\left( T^2(\mu)(T(Z)(X(m)))\subtraction\limits_{T(p_A)} 
T(Z^{\hlift})(T(\mu)(X(m))) \right).
\end{equation*}
Focus on $T^2(\mu)(T(Z)(X(m)))\subtraction\limits_{T(p_A)}
T(Z^{\hlift})(T(\mu)(X(m)))$. We can rewrite this as
\begin{equation}\label{second_line}
T (T(\mu)\circ Z )(X(m)) \subtraction\limits_{T(p_A)}
T(Z^{\hlift}\circ \mu )(X(m)). 
\end{equation}
At this point, we use Proposition~\ref{prop:Tw} below, that the warp of the tangent
of a grid is the tangent of the warp. 
We apply this to the grid on the Down face which is as given in (\ref{eq:GT}). 
From Example~\ref{ex:conn} 
the warp of (\ref{eq:GT}) is $\nabla_Z(\mu)$ (see (\ref{eq:conn})).
The tangent of (\ref{eq:GT}) is
\begin{equation*}
\begin{tikzcd}[row sep=2cm, column sep = 2cm]
T^2A \arrow[r,"T^2(q)",swap]\arrow[r,<-, bend left = 25,
"T^2(\mu)"]\arrow[d,"T(p_A)"]\arrow[d,<-, bend right = 25,swap, "T(Z^{\hlift})"]
&T^2M \arrow[d,"T(p)",swap]\arrow[d,<-, bend left = 25, "T(Z)"]\\
TA\arrow[r,"T(q)"]\arrow[r,<-, bend right = 25,swap, "T(\mu)"]
&TM,
\end{tikzcd}
\end{equation*}
and so its warp is given, for any $x\in T_mM$, by Proposition~\ref{prop:Tw},
\begin{equation}\label{eqn:tangent-grid}
(T^2(\mu)\circ T(Z))(x)\subtraction\limits_{T(p_A)}
(T(Z^{\hlift})\circ T(\mu))(x) =
\overline{T(\nabla_Z\mu)(x)}^{U}\addition\limits_{T^2(q)}
T(\wtilde{0}^{TA})(T(\mu)(x)).
\end{equation}
Here we have denoted by $\wtilde{0}^{TA}$ the zero section of
$TA\xrightarrow{p_A} A$. For $x = X(m)$, the right hand side of (\ref{eqn:tangent-grid}) is
equal to (\ref{second_line}). Therefore, (\ref{second_line}) is equal to
\begin{equation}\label{third_line}
\overline{T(\nabla_Z\mu)(X(m))}^{U}\addition\limits_{T^2(q)}
T(\wtilde{0}^{TA})(T(\mu)(X(m))).
\end{equation}

We return now to our calculation of (\ref{first_line}). Applying
$J_A$ to (\ref{third_line}), we have that (\ref{first_line}) is
\begin{equation*}
J_A\left(\overline{T(\nabla_Z\mu)(X(m))}^{U}\right)\addition\limits_{T^2(q)} 
J_A\left(T(\wtilde{0}^{TA})(T(\mu)(X(m)))\right). 
\end{equation*}
The addition over $T^2(q)$ does not change under $J_A$, by Lemma
\ref{lem:Jadd}. From (\ref{J_A:(U-D)to(B-F)}) we have
\begin{equation*}
J_A\left(\overline{T(\nabla_Z\mu)(X(m))}^{U}\right)\addition\limits_{T^2(q)} 
J_A\left(T(\wtilde{0}^{TA})(T(\mu)(X(m)))\right) =
\overline{T(\nabla_Z\mu)(X(m))}^{B} \addition\limits_{T^2(q)}
0^{T^2A}_{T(\mu)(X(m))},
\end{equation*}
This completes the calculation of the warp of the Back face; taking into 
consideration the orientation of the Back face, the warp is $-T(\nabla_Z\mu)\in
\Gamma_{TM}TA$. 

\begin{proposition}
\label{prop:Tw}
Let $(\xi,X)$ and $(\eta,Y)$ be a grid on a \dvb $D$ with warp
$\warp(\xi,\eta)\in\Ga C$. Then $(T(\xi),T(X))$ and $(T(\eta),T(Y))$ form a grid 
on the \dvb $TD$ in {\rm (\ref{tangent_grid_of_TD})} below and the warp of the 
tangent grid $(T(\xi),T(X))$, $(T(\eta),T(Y))$ is
$T(\warp(\xi,\eta))\in\Ga_{TM}(TC)$.
\end{proposition}

\begin{equation}\label{tangent_grid_of_TD}
\begin{tikzcd}[row sep=1.5cm, column sep = 1.5cm]
D \arrow[r]\arrow[r,<-, bend left = 25, "\xi"]\arrow[d]\arrow[d,<-,
bend right = 25,swap, "\eta"] &B \arrow[d]\arrow[d,<-, bend left =
25, "Y"]\\
A\arrow[r]\arrow[r,<-, bend right = 25,swap, "X"]
&M,
\end{tikzcd}
\qquad\qquad
\begin{tikzcd}[row sep=1.5cm, column sep = 1.5cm]
TD \arrow[r]\arrow[r,<-, bend left = 25, "T(\xi)"]\arrow[d]\arrow[d,<-,
bend right = 25,swap, "T(\eta)"] &TB \arrow[d]\arrow[d,<-, bend left =
25, "T(Y)"]\\
TA\arrow[r]\arrow[r,<-, bend right = 25,swap, "T(X)"]
&TM.
\end{tikzcd}
\end{equation}

\begin{proof}
We leave the verification that $TD$ is a \dvb and that $(T(\xi),T(X))$ and 
$(T(\eta),T(Y))$ are linear sections of $TD$ to the reader. 

The warp $\warp(\xi,\eta)\in \Gamma C$ of $(\xi,X)$, $(\eta,Y)$ 
is given as usual by
\begin{equation*}
(\xi\circ Y)(m)\subtraction\limits_{A}(\eta\circ X)(m) =
0^D_{X(m)}\addition\limits_B \warp(\xi,\eta)(m),\quad
(\xi\circ Y)(m)\subtraction\limits_{B}(\eta\circ X)(m) =
0^D_{Y(m)}\addition\limits_A \warp(\xi,\eta)(m).
\end{equation*}

We calculate the warp of the tangent grid. From the definition of a warp, for
$x\in T_mM$,
\begin{equation*}
(T(\xi)\circ T(Y))(x)\subtraction\limits_{TA}
(T(\eta)\circ T(X))(x) = 
T(0^D\circ X)(x)\addition\limits_{TB}\warp(T(\xi), T(\eta))(x).
\end{equation*}
Write $x = \frac{d}{dt}m_t\Big\rvert_{t = 0}$, for $m_t$ a curve in $M$
with tangent vector $x$ at $t = 0$.
Then, for $F\in C^{\infty}(D)$,
\begin{eqnarray}
\label{eq:Twcalc}
\left(T(\xi)\circ T(Y))(x)\subtraction\limits_{TA}
(T(\eta)\circ T(X)\right)(x)(F) &
 = &
\frac{d}{dt}F\left((\xi\circ Y)(m_t)\subtraction\limits_A (\eta\circ X)(m_t) 
\right)\Big\rvert_{t = 0}\nonumber\\
& = & \frac{d}{dt}F\left(0^D_{X(m_t)}\addition\limits_B \warp(\xi,\eta)(m_t)
\right)\Big\rvert_{t  = 0}\nonumber\\
& = & \frac{d}{dt}F\left((0^D\circ X)(m_t)\addition\limits_B
\warp(\xi,\eta)(m_t) \right)\Big\rvert_{t  = 0}\nonumber\\
& = & \left(T(0^D\circ X)(x)\addition\limits_{TB} T(\warp(\xi,\eta))(x)
\right)(F). 
\end{eqnarray}
Here we used the formula
\begin{equation}
\label{addition_over_TA}
(T(\xi)(\Phi_1)\addition\limits_{TA}T(\xi)(\Phi_2))(F) = 
\frac{d}{dt}F(\xi\circ \ph^1_t\addition\limits_A \xi\circ \ph^2_t)\Big\rvert_{t = 0}
\end{equation}
for the addition in $TD\to TA$, and the corresponding formula for scalar multiplication. 
In (\ref{addition_over_TA}) $F\in C^{\infty}(D)$, the $\Phi_i$ are elements of $TB$ and 
the $\ph^i_t$ are curves in $B$ with 
$$
\Phi_i = \frac{d}{dt}\ph^i_t\Big\rvert_{t = 0},\quad i = 1,2.
$$
Given $T(q_B)(\Phi_1) = T(q_B)(\Phi_2)$ we can arrange that 
$q_B(\ph^1_t) = q_B(\ph^2_t)$ for $t$ near zero. 

By uniqueness of the core element, it follows from (\ref{eq:Twcalc}) that 
\begin{equation*}
\warp(T(\xi), T(\eta))(x) = T(\warp(\xi,\eta))(x).
\end{equation*}
\end{proof}

\subsection{The three ultrawarps}

We now focus on the grids defined on the core \dvbs of $T^2A$.
We present a table with the results here, and outline the
calculations in the following subsections.

\medskip

\begin{table}[h]
\begin{tabular}{|c|c|c|c|}
\hline
Back & Front &  (B-F) core \dvb & $\uwarp_\baf$\\
\hline
&&&\\
$-T(\nabla_Z\mu)$ & $-\nabla_Z\mu$ & \begin{tikzcd}[row sep=1.5cm, column sep =
1.5cm] TA \arrow[r]\arrow[r,<-, bend left = 25,
"X^{\hlift}"]\arrow[d]\arrow[d,<-, bend right = 25,swap, "-T(\nabla_Z\mu)"] &A
\arrow[d]\arrow[d,<-, bend left = 25, "-\nabla_Z(\mu)"]\\
TM\arrow[r]\arrow[r,<-, bend right = 25,swap, "X"]
&M
\end{tikzcd} & $-\nabla_X\nabla_Z\mu$\\ 
&&&\\
\hline 
Left & Right & (L-R) core \dvb & $\uwarp_\lr$\\
\hline
&&&\\
$[Z^{\hlift},X^{\hlift}]$ & $[Z,X]$ & \begin{tikzcd}[row sep=1.5cm, column sep =
1.5cm] TA \arrow[r]\arrow[r,<-, bend left = 25,
"T(\mu)"]\arrow[d]\arrow[d,<-, bend right = 25,swap, "\leftbracket
Z^{\hlift}\comma X^{\hlift}\rightbracket"] &TM \arrow[d]\arrow[d,<-, bend left =
25, "\leftbracket Z\comma X\rightbracket"]\\
A\arrow[r]\arrow[r,<-, bend right = 25,swap, "\mu"]
&M
\end{tikzcd} & $-\nabla_{[Z,X]}\mu+R_{\nabla}(Z,X)(\mu)$\\
&&&\\
\hline
Up & Down & (U-D) core \dvb & $\uwarp_\ud$\\
\hline 
&&&\\
$T(\nabla_X\mu)$ &$\nabla_X\mu$ &
\begin{tikzcd}[row sep=1.5cm, column sep = 1.5cm]
TA \arrow[r]\arrow[r,<-, bend left = 25, "Z^{\hlift}"]\arrow[d]\arrow[d,<-,
bend right = 25,swap, "T(\nabla_X\mu)"] &A \arrow[d]\arrow[d,<-, bend left =
25, "\nabla_X\mu"]\\
TM\arrow[r]\arrow[r,<-, bend right = 25,swap, "Z"]
&M
\end{tikzcd} & $\nabla_Z\nabla_X\mu$\\
\hline
\end{tabular}
\caption{\label{tbl:T2A}}
\end{table}

\subsubsection{Back-Front}

The Back face is the tangent \dvb of the prolonged bundle 
$TA\rightarrow TM$ and by the results of subsection~\ref{subsect:wbf}
we obtain
\begin{equation}
\label{eq:bfw}
T^2(\mu)\circ\widetilde{Z} - \wtilde{Z^{\hlift}}^A\circ T(\mu)\,\ced{}\,T(\nabla_Z\mu).
\end{equation}

Taking into account the orientation of the Back face, the warp is 
$-T(\nabla_Z\mu)\in\Ga_{TM}(TA)$. For the Front face, with the appropriate
orientation, the warp is $-\nabla_Z\mu\in\Ga A$, by Example~\ref{ex:conn}. 
Therefore the ultrawarp for the Back-Front core \dvb 
(first row of Table~\ref{tbl:T2A}) is, again using Example~\ref{ex:conn} 
\begin{equation*}
-T(\nabla_Z\mu)\circ X +
X^{\hlift}(\nabla_Z\mu)\,\ced{}\, -\nabla_X\nabla_Z\mu.
\end{equation*}

\subsubsection{Left-Right}
The Left face is the tangent double vector bundle $T^2A$ for the manifold $A$. 
We therefore apply (\ref{eq:AMRi}). Taking into account the orientation of the
Left face, we have
$$
T(X^{\hlift})\circ Z^{\hlift} - \wtilde{Z^{\hlift}}^A\circ X^{\hlift}\,\ced{}\,
[Z^{\hlift}, X^{\hlift}]. 
$$
The Right face is $T^2M$ so the warp is $[Z,X]\in\X(M)$. So the warp of the core
\dvb in the second row of Table~\ref{tbl:T2A} is defined by
\begin{equation}\label{LR-ultrawarp-T2A}
 T(\mu)\circ [Z,X] - [Z^\hlift,X^\hlift]\circ \mu.
\end{equation}
This grid is of a new type; we encounter it here for the first time. 

First, what is $[Z^{\hlift},X^{\hlift}]$? When a connection in a vector bundle is
formulated in terms of a lifting of vector fields from the base to the total space
(see Example~\ref{ex:conn}), the difference $[Z^{\hlift},X^{\hlift}] - [Z,X]^{\hlift}$ 
is related to 
the curvature of the connection. This is analogous to 
describing the curvature of a connection in a principal bundle $P(M,G)$ as the failure 
of the horizontal distribution to be involutive. In \cite{Dieudonne:III} and other works 
of that period, it is stated that this expression reduces to the endomorphism $R_\nabla(X,Y)$. 
We deduce this now from the warp theorem. 

Both $[Z^{\hlift},X^{\hlift}]$ and $[Z,X]^{\hlift}$ project to $[Z,X]$ and therefore
their difference is a linear and vertical vector field on $A$.

At this point we take a closer look at linear sections $(\eta, 0^B)$,
where $\eta\in \Gamma_A D$ and $0^B\in\Ga B$, of a general \dvb $D$. It can
be shown that such a section $\eta$ induces a \vb map $\ph: A\rightarrow C$,
over $M$, and for $a\in A$,
\begin{equation*}
\eta(a) = \ph(a)\addition\limits_{B}0^D_a.
\end{equation*}
We denote such linear sections by $(\ph^{\bolt},0^B)$, where $\ph$ is the
induced \vb map.
It follows immediately that the warp of a horizontal linear section $(\xi,X)$
and of a vertical linear section over the zero section $(\ph^{\bolt}, 0^B)$ is
\begin{equation}\label{warp-of-linear-and-bolt}
\warp(\xi,\ph^{\bolt}) = -\ph\circ X.
\end{equation}
Back to (\ref{LR-ultrawarp-T2A}). We will need the following Proposition.
\begin{proposition}\label{warp-of-sum=sum-of-warp}
Given two grids $(\xi,X)$, $(\eta,Y)$, and $(\xi,X)$, $(\ph^{\bolt},0^B)$ on
a \dvb $D$, 
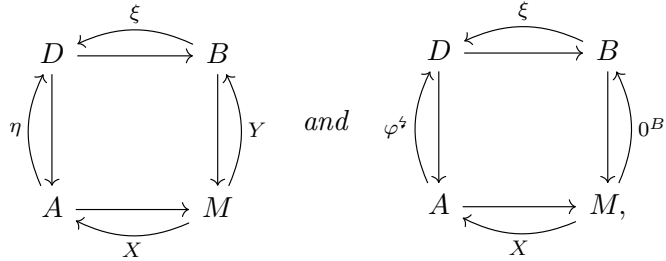
\begin{equation*}
\begin{tikzcd}[row sep=1.5cm, column sep = 1.5cm]
D \arrow[r]\arrow[r,<-, bend left = 25, "\xi"]\arrow[d]\arrow[d,<-,
bend right =25,swap, "\eta"] &B \arrow[d]\arrow[d,<-, bend left =
25, "Y"]\\
A\arrow[r]\arrow[r,<-, bend right = 25,swap, "X"]
&M
\end{tikzcd}\quad\text{and}\quad
\begin{tikzcd}[row sep=1.5cm, column sep = 1.5cm]
D \arrow[r]\arrow[r,<-, bend left = 25, "\xi"]\arrow[d]\arrow[d,<-,
bend right = 25,swap, "\ph^{\bolt}"] &B \arrow[d]\arrow[d,<-, bend left =
25, "0^B"]\\
A\arrow[r]\arrow[r,<-, bend right = 25,swap, "X"]
&M,
\end{tikzcd}
\end{equation*}
with warps $\warp(\xi,\eta)$ and $\warp(\xi,\ph^{\bolt})$, then
\begin{equation*}
\warp(\xi,\eta\addition\limits_A\ph^{\bolt}) =
\warp(\xi,\eta)+\warp(\xi,\ph^{\bolt}).
\end{equation*}
\end{proposition}
We leave the proof to the reader.

We can now rewrite the grid on the (L-R) core \dvb as the sum of 
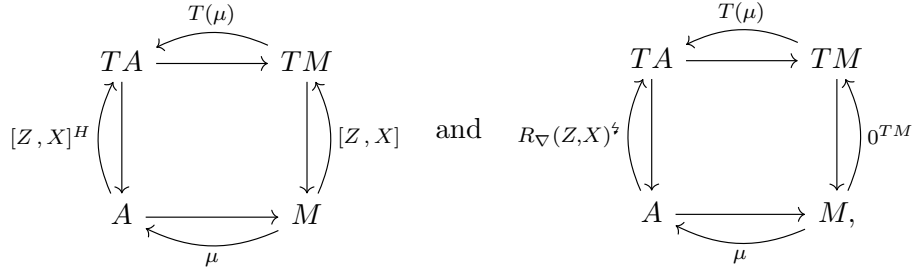
\begin{equation*}
\begin{tikzcd}[row sep=1.5cm, column sep = 1.5cm]
TA \arrow[r]\arrow[r,<-, bend left = 25, "T(\mu)"]\arrow[d]\arrow[d,<-,
bend right = 25,swap, "\leftbracket Z\comma X\rightbracket^{\hlift}"]
&TM \arrow[d]\arrow[d,<-, bend left = 25, "\leftbracket
Z\comma X\rightbracket"]\\
A\arrow[r]\arrow[r,<-, bend right = 25,swap, "\mu"]
&M
\end{tikzcd}\quad\text{and}\quad
\begin{tikzcd}[row sep=1.5cm, column sep = 1.5cm]
TA \arrow[r]\arrow[r,<-, bend left = 25, "T(\mu)"]\arrow[d]\arrow[d,<-,
bend right = 25,swap, "{R_{\nabla}(Z,X)}^{\bolt}"]
&TM \arrow[d]\arrow[d,<-, bend left = 25, "0^{TM}"]\\
A\arrow[r]\arrow[r,<-, bend right = 25,swap, "\mu"]
&M,
\end{tikzcd}
\end{equation*}
so (\ref{LR-ultrawarp-T2A}) is now, from Proposition
\ref{warp-of-sum=sum-of-warp},
\begin{equation*}
\left(T(\mu)\circ [Z,X]-[Z,X]^{\hlift}\circ \mu\right)+
\left(T(\mu)\circ 0^{TM} - (R_{\nabla}(Z,X)^{\bolt})\circ \mu\right).
\end{equation*}
From Example~\ref{ex:conn} 
\begin{equation*}
T(\mu)\circ [Z,X]-[Z,X]^{\hlift}\circ \mu \,\ced{}\, \nabla_{[Z,X]}\mu,
\end{equation*}
and from (\ref{warp-of-linear-and-bolt})
\begin{equation*}
T(\mu)\circ 0^{TM} - (R_{\nabla}(Z,X)^{\bolt})\circ \mu \,\ced{}\,
-R_{\nabla}(Z,X)(\mu).
\end{equation*}
So in total, the warp of this core \dvb will be 
\begin{equation*}
\nabla_{[Z,X]}\mu - R_{\nabla}(Z,X)(\mu).
\end{equation*}
Taking into consideration the orientation of the core \dvb,
take the opposite sign
\begin{equation*}
-\nabla_{[Z,X]}\mu + R_{\nabla}(Z,X)(\mu).
\end{equation*}

\subsubsection{Up-Down}

The warp of the Down face is, again by Example~\ref{ex:conn}, 
$\nabla_X\mu$. 
For the warp of the Up face, we use Proposition~\ref{prop:Tw}, and obtain
$T(\nabla_X\mu)\in\Ga_{TM}(TA)$.
Therefore the warp of the grid in the third row of Table~\ref{tbl:T2A} is
\begin{equation*}
\nabla_Z\nabla_X\mu.
\end{equation*}

This completes the exposition of Table~\ref{tbl:T2A}. 

\bigskip

The warp theorem now gives us that
\begin{equation}
-\nabla_X\nabla_Z\mu - \nabla_{[Z,X]}\mu + R_\nabla(Z,X)(\mu) +\nabla_Z\nabla_X\mu = 0.
\end{equation}
This is the definition of the curvature of $\nabla$ via differential
operators. Therefore, we see that if we start with the concept of a
connection $\nabla$, and apply the warp theorem to the grid in $T^2A$,
we are led to define the quantity $R_\nabla(Z,X)(\mu)$ in this way.

\section{The triple tangent bundle $T^3M$ and the Jacobi identity}
\label{sect:Jac}

In this section we consider the triple tangent bundle $T^3M$ of a manifold $M$
and construct a grid on it, for which the Jacobi identity emerges as a
consequence of the warp theorem. A version of this approach was given by
Mackenzie \cite{Mackenzie:pJihw}. We present here a clearer and more
detailed calculation. 

Take $E$ to be $T^3M$, the triple tangent bundle.
This is a special case of $T^2A$, for $A = TM$:
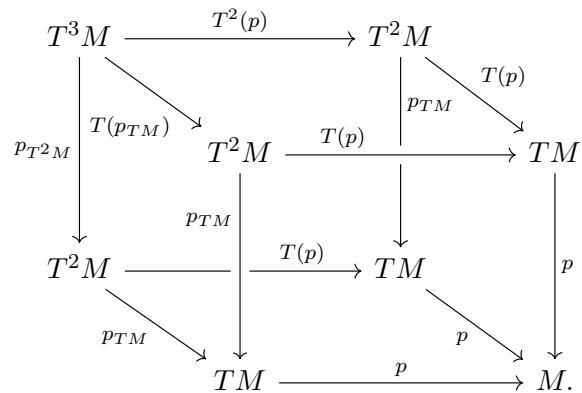
\begin{equation*}
\begin{tikzcd}[row sep=2.45em, column sep=2.45em]
T^3 M \arrow[rr,"T^2(p)"] \arrow[dr,swap,"T(p_{TM})",near end]
\arrow[dd,swap,"p_{T^2M}"] && T^2M \arrow[dd,"p_{TM}" near start] \arrow[dr,"T(p)"] \\
& T^2M \arrow[rr,crossing over,"T(p)" near start] &&
  TM \arrow[dd,"p"] \\
T^2M \arrow[rr,"T(p)" near end] \arrow[dr,swap,"p_{TM}"] && TM
\arrow[dr,swap,"p"]
\\
& TM \arrow[rr,"p"] \arrow[uu,<-,crossing over,"p_{TM}" near end]&& M.
\end{tikzcd}
\end{equation*} 
The three lower faces are copies of $T^2M$. The Left face is the double
tangent bundle of the manifold $TM$. The Back face is not a double tangent
bundle; it is the tangent \dvb of \newline $T^2M\xrightarrow{T(p)}TM$. 
The Up face is obtained by applying the tangent functor to $T^2M$.

Starting with three vector fields $X$, $Y$, and $Z$, each a section of one of
the three copies of $TM$, one can build a grid on $T^3M$ as follows; see (\ref{eq:T3Mgrid}) below. 
\begin{itemize}
  \item The front-back linear double section
  $(T(\wtilde{X});T(X),\wtilde{X};X)$. Take the complete lift of $X$
  across the Down face, and apply the tangent functor to the linear section $(\wtilde{X},X)$.
  \item The right-left linear double section $(T^2(Y);T(Y),T(Y);Y)$.
  Apply the tangent functor to $Y$ and then to $T(Y)$. 
  \item The down-up linear double section
  $(\wtilde{\wtilde{Z}};\wtilde{Z},\wtilde{Z};Z)$. Take the complete lift of $Z$
    across the Front face, and the complete lift of this across the Left face.
    Likewise take the complete lift of $Z$ across the Right face. One does need
    to check that $(\wtilde{\wtilde{Z}},\wtilde{Z})$ is indeed a linear section
    of the Back face.
\end{itemize}
Equation (\ref{eq:T3Mgrid}) shows the entire grid. 
\begin{equation}
  \label{eq:T3Mgrid}
\begin{tikzcd}[row sep=2.45em, column sep=2.45em]
T^3M \arrow[rr,<-,"T^2(Y)"] \arrow[dr,<-,swap,"T(\widetilde{X})",near end]
\arrow[dd,<-,swap,"\widetilde{\widetilde{Z}}"] && T^2M
\arrow[dd,<-,"\widetilde{Z}" near start] \arrow[dr,<-,"T(X)"] \\
& T^2M \arrow[rr,<-,crossing over,swap,"T(Y)" near start] &&
TM \arrow[dd,<-,"Z"] \\
T^2M \arrow[rr,<-,swap,"T(Y)" near start] \arrow[dr,<-,swap,"\widetilde{X}"] &&
TM\arrow[dr,<-,"X"]\\
& TM \arrow[rr,<-,swap,"Y"] \arrow[uu,crossing over,"\widetilde{Z}" near end]&&
M.
\end{tikzcd} 
\end{equation}
We now calculate the three ultrawarps defined by this grid. To do this,
we calculate the core morphisms of the three linear double sections, and
the warps of the six faces.

First, the core morphisms. These follow in an analogous way as in the example of
$T^2A$,
\begin{itemize}
  \item The core morphism of $(T(\wtilde{X});T(X),\wtilde{X};X)$
  is $(\wtilde{X},X)$.
  \item The core morphism of $(T^2(Y);T(Y),T(Y);Y)$
  is $(T(Y),Y)$.
  \item The core morphism of
  $(\wtilde{\wtilde{Z}};\wtilde{Z},\wtilde{Z};Z)$ is $(\wtilde{Z},Z)$.
\end{itemize}

To calculate the warps of the six faces, we take into consideration the
orientation of the faces of a \tvb. 
For the lower faces, by (\ref{eq:AMRi}):
\begin{itemize}
  \item For the Front face:\qquad
 $\widetilde{Z}(Y)\subtraction T(Y)(Z)\,\ced{}\,[Y,Z]$.  
\item For the Right face:\qquad
 $T(X)(Z)\subtraction \wtilde{Z}(X)\,\ced{}\,[Z,X]$.
\item For the Down face:\qquad
 $T(Y)(X)\subtraction \wtilde{X}(Y)\,\ced{}\,[X,Y]$.
\end{itemize}
We now calculate the warps of the upper faces.
\subsubsection*{Back face}
 The warp of the Back face, for $x\in TM$, is given by
\begin{equation}\label{warp-of-back-face}
\widetilde{\widetilde{Z}}\circ T(Y)(x)\subtraction\limits_{T^2(p)}T^2(Y)\circ
\widetilde{Z}(x) = \warp_{\text{back}}(x)\addition\limits_{p_{T^2M}}
\hat{0}_{\widetilde{Z}(x)}.
\end{equation}
As we mentioned, the Back face is the tangent \dvb of
$T^2M\xrightarrow{T(p)}TM$. Apply $T(J)$ to it, the tangent of the canonical
involution $J:T^2M\rightarrow T^2M$. The resulting \dvb is now the double
tangent bundle of $TM$.
In fact, $T(J)$ is a \tvb morphism, and maps the Back face of $T^3M$ to the
double tangent bundle of $TM$ as shown in (\ref{T(J) maps Back to Left face}). 
\begin{equation}\label{T(J) maps Back to Left face}
\begin{tikzcd}[row sep=3em, column sep=3em]
T^3M\arrow[rr,"T^2(p)"] \arrow[red,dr,"T(J)"]
\arrow[dd,swap,"p_{T^2M}"] && T^2M \arrow[dd,"p_{TM}" near end]
\arrow[red,dr,"\id"] \\
& T^3M \arrow[rr,crossing over,"T(p_{TM})" near start] &&
  T^2M \arrow[dd,"p_{TM}"] \\
T^2M \arrow[rr,"T(p)" near end] \arrow[red,dr,swap,"J"] && TM
\arrow[red,dr,swap,"\id"]
\\
& T^2M \arrow[rr,"p_{TM}"] \arrow[uu,<-,crossing over,"p_{T^2M}" near end]&& TM.
\end{tikzcd}
\end{equation}
Also, the core morphism of (\ref{T(J) maps Back to Left face}) is $(J,\id)$. 
Hence, applying $T(J)$ to (\ref{warp-of-back-face}),
\begin{equation}\label{T(J) of warp of back face}
T(J)\left(\widetilde{\widetilde{Z}}\circ
T(Y)(x)\subtraction\limits_{T^2(p)}T^2(Y) \circ\widetilde{Z}(x)\right) =
 J(\warp_{\text{back}}(x))\addition\limits_{p_{T^2M}}\hat{0}_{\wtilde{Z}(x)}.
\end{equation}
Note that $T(J)$ changes the \vb structure
over which the subtraction of the left hand side takes place, and
$\subtraction\limits_{T^2(p)}$ will become $\subtraction\limits_{T(p_{TM})}$.
Applying $T(J)$ to the grid of the Back face yields the following grid on the
double tangent bundle of $TM$. 
\begin{equation*}
\begin{tikzcd}[row sep=2cm, column sep = 2cm]
T^3M\arrow[r,"T(p_{TM})",swap]\arrow[r,<-, bend left = 20,
"T(\wtilde{Y})"]\arrow[d,"p_{T^2M}"]\arrow[d,<-, bend right = 20,swap,
"\wtilde{\wtilde{Z}}"] &T^2M \arrow[d,"p_{TM}",swap]\arrow[d,<-, bend left = 20,
"\wtilde{Z}"]\\
T^2M\arrow[r,"p_{TM}"]\arrow[r,<-, bend right = 20,swap, "\wtilde{Y}"]
&TM.
\end{tikzcd}
\end{equation*}
Therefore, expanding the left hand side of (\ref{T(J) of warp of back face}),
 \begin{multline*}
 T(J)\left(\widetilde{\widetilde{Z}}\circ
 T(Y)(x)\subtraction\limits_{T^2(p)}T^2(Y) \circ\widetilde{Z}(x)\right) = 
 T(J) ((\widetilde{\widetilde{Z}}\circ
 T(Y))(x))\subtraction\limits_{T(p_{TM})}T(J)((T^2(Y) \circ\widetilde{Z})(x))
 \\= (\widetilde{\widetilde{Z}}\circ
 \wtilde{Y})(x)\subtraction\limits_{T(p_{TM})} (T(\wtilde{Y})\circ
 \wtilde{Z})(x) \equals\limits^{(\ref{eq:AMRi})}
 -[\wtilde{Z},\wtilde{Y}](x)\addition\limits_{p_{T^2M}}\hat{0}_{\wtilde{Z}(x)} = 
 \wtilde{[Y,Z]}(x)\addition\limits_{p_{T^2M}}\hat{0}_{\wtilde{Z}(x)}.
 \end{multline*}
Substituting this into (\ref{T(J) of warp of back face}),
\begin{equation*}
\wtilde{[Y,Z]}(x)\addition\limits_{p_{T^2M}}\hat{0}_{\wtilde{Z}(x)} = 
J(\warp_{\text{back}}(x))\addition\limits_{p_{T^2M}}\hat{0}_{\wtilde{Z}(x)},
\end{equation*}
and using that $J^2 = \id$, we obtain
\begin{equation*}
\warp_{\text{back}}= T([Y,Z]).
\end{equation*}
\subsection*{Left face}
The Left face is the double tangent bundle of $TM$, so we simply
  apply (\ref{eq:AMRi}) for the grid
  $(T(\wtilde{X}),\wtilde{X})$, $(\wtilde{\wtilde{Z}},\wtilde{Z})$,
  \begin{equation*}
T(\widetilde{X})\circ
\widetilde{Z}-\widetilde{\widetilde{Z}}\circ
\widetilde{X}\,\ced{}\,[\widetilde{Z}, \widetilde{X}] = \widetilde{[Z,X]},
\end{equation*}
  so $\warp_{\text{left}} = \wtilde{[Z,X]}$.
\subsubsection*{Up face}
For the Up face, using Proposition~\ref{prop:Tw}, it follows directly that
$\warp_{\text{up}} = T([X,Y])$.

\subsubsection{The three ultrawarps}
The three core \dvbs are all copies of $T^2M$, and their ultracore is
$TM\rightarrow M$.

The three core \dvbs in the usual order (B-F), (L-R),
and (U-D), with the induced grids from the original grid on $T^3M$,
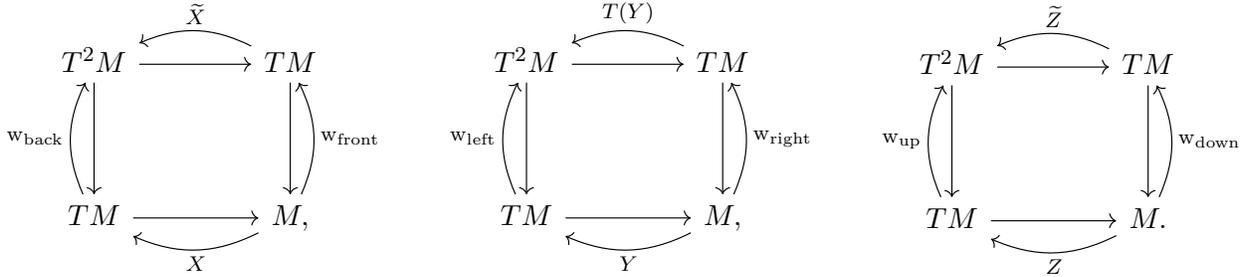
\begin{equation*}
\begin{tikzcd}[row sep=1.5cm, column sep = 1.5cm]
T^2M \arrow[r]\arrow[r,<-, bend left = 25, "\widetilde{X}"]\arrow[d]\arrow[d,<-,
bend right = 25,swap, "\warp_{\text{back}}"] &TM \arrow[d]\arrow[d,<-, bend left
= 25, "\warp_{\text{front}}"]\\
TM\arrow[r]\arrow[r,<-, bend right = 25,swap, "X"]
&M,
\end{tikzcd}\qquad
\begin{tikzcd}[row sep=1.5cm, column sep = 1.5cm]
T^2M \arrow[r]\arrow[r,<-, bend left = 25, "T(Y)"]\arrow[d]\arrow[d,<-,
bend right = 25,swap, "\warp_{\text{left}}"] &TM \arrow[d]\arrow[d,<-, bend left
= 25, "\warp_{\text{right}}"]\\
TM\arrow[r]\arrow[r,<-, bend right = 25,swap, "Y"]
&M,
\end{tikzcd}\qquad
\begin{tikzcd}[row sep=1.5cm, column sep = 1.5cm]
T^2M \arrow[r]\arrow[r,<-, bend left = 25, "\widetilde{Z}"]\arrow[d]\arrow[d,<-,
bend right = 25,swap, "\warp_{\text{up}}"] &TM \arrow[d]\arrow[d,<-, bend left
= 25, "\warp_{\text{down}}"]\\
TM\arrow[r]\arrow[r,<-, bend right = 25,swap, "Z"]
&M.
\end{tikzcd}
\end{equation*}
Finally, by (\ref{eq:AMRi}), the ultracore elements are
\begin{eqnarray*}
\warp_{\text{back}}\circ X-\widetilde{X}\circ
\warp_{\text{front}} =& T([Y,Z])\circ X-\widetilde{X}\circ [Y,Z]&\,\ced{}\,
[X,[Y,Z]],\\
\warp_{\text{left}}\circ Y - T(Y)\circ
\warp_{\text{right}} = &\widetilde{[Z,X]}\circ Y - T(Y)\circ [Z,X]& \,\ced{}\,
-[[Z,X],Y] = [Y,[Z,X]],\\
\warp_{\text{up}}\circ Z - \widetilde{Z}\circ
\warp_{\text{down}} = &T([X,Y])\circ Z-\widetilde{Z}\circ [X,Y]&\,\ced{}\,
[Z,[X,Y]].
\end{eqnarray*}
We see that in this way we have formulated the three terms of the Jacobi identity.
And applying the warp theorem, we obtain a conceptual proof of the 
Jacobi identity.

\newpage

\section{Appendix: warps and duality}
\label{sect:app}

We include here an alternative formula for the warp which relies on the duality
of \dvbs. We will need only a few basics from duality theory; for further
details see \cite[\S9.2]{Mackenzie:GT}.

Start with a \dvb $D$ and take the usual dual of $D$ with respect to $A$. This
defines another \dvb $D\duer A$, with core $B^*$. In a similar way we can define
the \dvb $D\duer B$, with core $A^*$. These duals are shown in (\ref{eq:duals}). 
\begin{equation}
\label{eq:duals}
\begin{tikzcd}[row sep=2cm, column sep = 2cm]
D \arrow[r, "q^D_B"] \arrow[d,"q^D_A"] & B \arrow[d,"q_B"] \\
A\arrow[r,"q_A"] &M,
\end{tikzcd}
\hspace*{12mm}
\begin{tikzcd}[row sep=2cm, column sep = 2cm]
D\duer A \arrow[r, "\gamma^A_{C^*}"] \arrow[d,"\gamma^A_A"] &C^* \arrow[d,""] \\
A\arrow[r,""] &M,
\end{tikzcd}
\hspace*{12mm}
\begin{tikzcd}[row sep=2cm, column sep = 2cm]
D\duer B \arrow[r, "\gamma^B_B"] \arrow[d,"\gamma^B_{C^*}"] & B \arrow[d,""] \\
C^* \arrow[r,""] &M.
\end{tikzcd}
\end{equation}
There exists a nondegenerate pairing between $D\duer A$ and $D\duer B$ over
$C^*$ \cite[9.2.2]{Mackenzie:GT}, which is natural up to sign. 
For $\Phi\in D\duer A$, and $\Psi\in D\duer
B$, with outlines $(\Phi;a,\kappa;m)$ and $(\Psi;\kappa,b;m)$ respectively, and for
any $d\in D$ with outline $(d;a,b;m)$,
\begin{equation}\label{mackenzie_pairing_between_duals}
\|\Phi,\Psi\|_{C^*} = \langle \Phi,d\rangle_A-\langle \Psi, d\rangle_B.
\end{equation}
This pairing induces two \dvb isomorphisms, namely, 
\begin{equation}
\label{the_iso_Z_A}
Z_A: D\duer A\rightarrow D\duer B\duer C^*,\qquad 
\langle Z_A(\Phi), \Psi\rangle_{C^*} = \|\Phi,\Psi\|_{C^*},
\end{equation} 
and 
\begin{equation}
\label{the_iso_Z_B}
Z_B: D\duer B\rightarrow D\duer A\duer C^*,\qquad 
\langle Z_B(\Psi), \Phi\rangle_{C^*} = \|\Phi,\Psi\|_{C^*}.
\end{equation}

The following material on linear sections and duality can be found in more 
detail in \cite{Mackenzie:2011}.

A linear section $(\eta,Y)$ of $D$ as in (\ref{eq:dgrid})
induces a linear map $\ell_\eta\co D\duer A\to\R$. This map is automatically
linear with respect to $D\duer A\to A$ but is also linear with respect to the
other \vb structure, $D\duer A\to C^*$. It therefore induces a section of
$D\duer A\duer C^*\to C^*$ and this is linear over $Y$. Using
(\ref{mackenzie_pairing_between_duals}), (\ref{the_iso_Z_A}), and
(\ref{the_iso_Z_B}),
we will define
a new pairing between sections of $D\duer A\duer C^*\to C^*$ and sections of
$D\duer B\duer C^*\to C^*$. In the case where we begin with a grid $(\eta,Y)$
and $(\xi,X)$ on $D$ as usual, we will show that this new pairing provides us
with a different way of expressing the warp $\warp(\xi,\eta)$.

How does a linear section $(\eta,Y)$ of $D$ define a linear section
$(\eta^{\sqcap},Y)$ of $D\duer A\duer C^*$?
Given $\kappa\in C^*$ we define 
$\eta^{\sqcap}(\kappa)\in D\duer A\duer C^*\Big\rvert_{\kappa}$ by defining
its pairing with $\Phi\in D\duer A\Big\rvert_{\kappa}$ to be
\begin{equation}\label{definition_of_eta_square_cap}
\langle \eta^{\sqcap}(\kappa),\Phi\rangle_{C^*} :=\ell_{\eta}(\Phi) =
\langle\Phi, \eta(\gamma^A_A(\Phi))\rangle_A.
\end{equation}
(We use notations such as $\Big\rvert_{\kappa}$ on \dvbs when the symbol
for the base point makes clear which structure is meant.) 

That $(\eta^{\sqcap}, Y)$ is a linear section of $D\duer A\duer C^*\rightarrow C^*$
follows immediately. The corresponding linear function is
$$
\ell_{\eta^{\sqcap}}: D\duer A \rightarrow  \R,\qquad
\Phi \mapsto  \langle \eta^{\sqcap}(\gamma^A_{C^*}(\Phi)),\Phi\rangle_{C^*},
$$
and of course $\ell_{\eta^{\sqcap}} = \ell_{\eta}$.
The following proposition states that $\ell_{\eta} = \ell_{\eta^{\sqcap}}$
is also linear with respect to $C^*$; see \cite[3.1]{Mackenzie:2011}.
\begin{proposition}
  \label{prop:Crelle}
If $(\eta,Y)$ is a linear section, then $\ell_{\eta}: D\duer A\rightarrow \R$ defined
by
\begin{equation*}
\Phi \mapsto \langle \Phi, \eta(\gamma^A_A(\Phi))\rangle
\end{equation*}
is linear with respect to $C^*$ as well as $A$, and the restriction of
$\ell_{\eta}$ to the core of $D\duer A$ is $\ell_Y:B^*\rightarrow \R$.
\end{proposition}

Therefore, we see that there exists a one-to-one correspondence between linear
sections $(\eta, X)$ of $D\rightarrow A$ and linear sections $(\eta^{\sqcap},
X)$ of $D\duer A\duer C^*\rightarrow C^*$.

Similarly, given $(\xi,X)$, define a section $\xi^{\sqcap}$ of $D\duer B\duer C^*\to C^*$ by
defining the pairing of $\xi^{\sqcap}(\kappa)$, for $\kappa\in C^*$, with 
$\Psi\in D\duer B\Big\rvert_{\kappa}$ to be
\begin{equation}\label{definition_of_xi_square_cap}
\langle \xi^{\sqcap}(\kappa),\Psi\rangle_{C^*} :=\ell_{\xi}(\Psi) =
\langle\Psi, \xi(\gamma^B_B(\Psi))\rangle_B.
\end{equation}
Now we want to define a pairing between $\xi^{\sqcap}\in \Gamma_{C^*}(D\duer B\duer C^*)$
and $\eta^{\sqcap}\in \Gamma_{C^*}(D\duer A\duer C^*)$:
\begin{equation*}
\lpair \xi^{\sqcap},\eta^{\sqcap}\rpair(\kappa) = \lpair 
\xi^{\sqcap}(\kappa),\eta^{\sqcap}(\kappa)\rpair
\end{equation*}
and it will follow that $\lpair \xi^{\sqcap},\eta^{\sqcap}\rpair\in
C^{\infty}(C^*)$.
Since $D\duer A\duer C^*$ and $D\duer A$ are dual \vbs over $C^*$, we have the usual
pairing between them. We will use it to define a pairing between
$D\duer A\duer C^*$ and $D\duer B\duer C^*$ over $C^*$. Using the map $Z_A^{-1}:
D\duer B\duer C^*\rightarrow D\duer A$ we define
\begin{equation}\label{pairing_via_Z_A}
\lpair \xi^{\sqcap}(\kappa),\eta^{\sqcap}(\kappa)\rpair = 
\langle \eta^{\sqcap}(\kappa),Z_A^{-1}(\xi^{\sqcap}(\kappa))\rangle_{C^*}.
\end{equation}

The following outlines are useful to help us keep track of the various
calculations,
\begin{eqnarray*}
\begin{tikzcd}[row sep=1.5cm, column sep = 1.5cm]
D\duer A\ni Z_A^{-1}(\xi^{\sqcap}(\kappa)) \arrow[r,mapsto," "]\arrow[d,mapsto,
swap," "] &\kappa \arrow[d,mapsto," "]\\
-X(m)\arrow[r,mapsto,swap," "]
&m, 
\end{tikzcd}&\quad&
\begin{tikzcd}[row sep=1.5cm, column sep = 1.5cm]
D\duer B\ni Z_B^{-1}(\eta^{\sqcap}(\kappa)) \arrow[r,mapsto," "]\arrow[d,mapsto,
swap," "] &Y(m)\arrow[d,mapsto," "]\\
\kappa\arrow[r,mapsto,swap," "]
&m, 
\end{tikzcd}\\
\begin{tikzcd}[row sep=1.5cm, column sep = 1.5cm]
\subtraction\limits_B \xi(Y(m))\arrow[r,mapsto," "]\arrow[d,mapsto,
swap," "] &Y(m) \arrow[d,mapsto," "]\\
-X(m)\arrow[r,mapsto,swap," "]
&m, 
\end{tikzcd}&\quad&
\begin{tikzcd}[row sep=1.5cm, column sep = 1.5cm]
\subtraction\limits_B \eta(X(m))\arrow[r,mapsto," "]\arrow[d,mapsto,
swap," "] &Y(m) \arrow[d,mapsto," "]\\
-X(m)\arrow[r,mapsto,swap," "]
&m.
\end{tikzcd}
\end{eqnarray*}
Note that the minus sign on $-X(m)$ of $Z_A^{-1}(\xi^{\sqcap}(\kappa))$ comes
from the fact that $Z_A$ induces $-\id_A:A\rightarrow A$ over $M$.

We can now begin calculations. Start with (\ref{pairing_via_Z_A}),
\begin{eqnarray*}
\lpair \xi^{\sqcap}(\kappa),\eta^{\sqcap}(\kappa)\rpair &=& 
\langle \eta^{\sqcap}(\kappa),Z_A^{-1}(\xi^{\sqcap}(\kappa))\rangle_{C^*}\\
&\equals\limits^{(\ref{definition_of_eta_square_cap})}&
\ell_{\eta}(Z_A^{-1}(\xi^{\sqcap}(\kappa))) \\
&= & -\langle Z_A^{-1}(\xi^{\sqcap}(\kappa)),\eta(X(m))\rangle_A.
\end{eqnarray*}
Now use (\ref{mackenzie_pairing_between_duals}), with $\Phi =
Z_A^{-1}(\xi^{\sqcap}(\kappa))$, $\Psi = Z_B^{-1}(\eta^{\sqcap}(\kappa))$, and 
$d = \subtraction\limits_B \eta(X(m))$,
\begin{multline*}
\|Z_A^{-1}(\xi^{\sqcap}(\kappa)),Z_B^{-1}(\eta^{\sqcap}(\kappa))\|_{C^*}
= \langle
Z_A^{-1}(\xi^{\sqcap}(\kappa)),\subtraction\limits_B\eta(X(m))\rangle_A -
\langle
Z_B^{-1}(\eta^{\sqcap}(\kappa)),\subtraction\limits_B\eta(X(m))\rangle_B\\\Rightarrow
-\langle Z_A^{-1}(\xi^{\sqcap}(\kappa)),\eta(X(m))\rangle_A =
\|Z_A^{-1}(\xi^{\sqcap}(\kappa)),Z_B^{-1}(\eta^{\sqcap}(\kappa))\|_{C^*} -
\langle Z_B^{-1}(\eta^{\sqcap}(\kappa)),\eta(X(m))\rangle_B
\end{multline*}
Returning to the previous calculations
\begin{eqnarray}\label{comparison}
\lpair \xi^{\sqcap}(\kappa),\eta^{\sqcap}(\kappa)\rpair &=& \nonumber
-\langle Z_A^{-1}(\xi^{\sqcap}(\kappa)),\eta(X(m))\rangle_A\\\nonumber
& = & \|Z_A^{-1}(\xi^{\sqcap}(\kappa)),Z_B^{-1}(\eta^{\sqcap}(\kappa))\|_{C^*} -
\langle Z_B^{-1}(\eta^{\sqcap}(\kappa)),\eta(X(m))\rangle_B\\
&\equals\limits^{(\ref{the_iso_Z_A})}&
\langle 
\xi^{\sqcap}(\kappa),Z_B^{-1}(\eta^{\sqcap}(\kappa))\rangle_{C^*}
-\langle Z_B^{-1}(\eta^{\sqcap}(\kappa)),\eta(X(m))\rangle_B\\
& \equals\limits^{(\ref{definition_of_xi_square_cap})} &\nonumber 
\langle Z_B^{-1}(\eta^{\sqcap}(\kappa)),\xi(Y(m))\rangle_B\nonumber
-\langle Z_B^{-1}(\eta^{\sqcap}(\kappa)),\eta(X(m))\rangle_B
\\\nonumber & = & \langle
Z_B^{-1}(\eta^{\sqcap}(\kappa)), \xi(Y(m))\subtraction\limits_B
\eta(X(m))\rangle_B \\\nonumber 
& = &  \langle Z_B^{-1}(\eta^{\sqcap}(\kappa)),\warp(\xi,\eta)(m)
\addition\limits_A 0^D_{Y(m)}\rangle_B.\nonumber
\end{eqnarray}
By equation (16) from \cite[p. 348]{Mackenzie:GT}, since
$\gamma^B_{C^*}(Z_B^{-1}(\eta^{\sqcap}(\kappa))) = \kappa$, we can rewrite the
last expression as
\begin{equation*}
\langle Z_B^{-1}(\eta^{\sqcap}(\kappa)),
\warp(\xi,\eta)(m)\addition\limits_A 0^D_{Y(m)}\rangle_B = \langle
\kappa,\warp(\xi,\eta)(m)\rangle = \ell_{\warp(\xi,\eta)}(\kappa).
\end{equation*}

In other words, we have shown that, for $\kappa\in C^*$, 
\begin{equation*}
\lpair\xi^{\sqcap},\eta^{\sqcap}\rpair(\kappa) = \ell_{\warp(\xi,\eta)}(\kappa). 
\end{equation*}
We state this result formally. 
\begin{proposition}
Let $(\xi, X)$ and $(\eta,Y)$ be linear sections forming a grid on a \dvb $D$. Then
\begin{equation*}
\lpair\xi^{\sqcap},\eta^{\sqcap}\rpair = \ell_{\warp(\xi,\eta)}. 
\end{equation*}
\end{proposition}

An equivalent way of defining the pairing (\ref{pairing_via_Z_A}) is via the map
$Z_B^{-1}: D\duer A\duer C^*\to D\duer B$, as follows
\begin{equation}\label{pairing_via_Z_B}
\lpair \xi^{\sqcap}(\kappa),\eta^{\sqcap}(\kappa)\rpair = 
\langle \xi^{\sqcap}(\kappa),Z_B^{-1}(\eta^{\sqcap}(\kappa))\rangle_{C^*}.
\end{equation}
It is easy to check that both (\ref{pairing_via_Z_A}) and
(\ref{pairing_via_Z_B}) define the same pairing. And comparing
(\ref{pairing_via_Z_B}) with (\ref{comparison}), we see that
\begin{equation}\label{term_zero}
\langle Z_B^{-1}(\eta^{\sqcap}(\kappa)),
\eta(X(m))\rangle_B = 0.
\end{equation}
Indeed, to see this, let $\Phi\in D\duer A$, with outline $(\Phi;
X(m),\kappa;m)$. Then, via (\ref{mackenzie_pairing_between_duals}), rewrite
(\ref{term_zero}) as follows
\begin{eqnarray*}
\langle Z_B^{-1}(\eta^{\sqcap}(\kappa)),
\eta(X(m))\rangle_B & = &
\langle \Phi, \eta(X(m))\rangle_A - \|\Phi,
Z_B^{-1}(\eta^{\sqcap}(\kappa))\|_{C^*}\\
& \equals\limits^{(\ref{the_iso_Z_B})}&
\langle \Phi, \eta(X(m))\rangle_A - 
\langle \eta^{\sqcap}(\kappa), \Phi\rangle_{C^*}\\
&\equals\limits^{(\ref{definition_of_eta_square_cap})}&
\langle \Phi, \eta(X(m))\rangle _A - \langle \Phi, \eta(X(m))\rangle _A\\
& = &0.
\end{eqnarray*}

\subsection{Example with $T^2M$}
In case of $(\wtilde{X}, X)$ and $(T(Y),Y)$, what are the corresponding
${\wtilde{X}}^{\sqcap}$ and $T(Y)^{\sqcap}$?
Starting with $(T(Y),Y)$,  
we will calculate $T(Y)^{\sqcap}$ using $\ell_{T(Y)}$,
$$ 
\ell_{T(Y)}:T^{\bullet}(TM) \to \R,\qquad
T^{\bullet}_xTM \ni \xi  \mapsto   \langle \xi, T(Y)(x)\rangle_{TM},
$$ 
where $x\in TM$ and $T^{\bullet}_xTM$ is the fibre of $T^{\bullet}(TM)$ over $x\in TM$.
Here $T^{\bullet}(TM)$ is the dual of $T(p): T^2M\to TM$ (see \cite[p. 355]{Mackenzie:GT} for
the notation). The function $\ell_{T(Y)}$ is linear with respect to both $TM$ and $T^*M$
(see Proposition~\ref{prop:Crelle}). Its linearity with respect to
$TM$ will lead us back to $T(Y)$. We are interested in its linearity with
respect to $T^*M$. This will define a linear section of the dual of the \vb
$T^{\bullet}(TM)\rightarrow T^*M$, that is, of $T^{\bullet}(TM)\duer
T^*M\rightarrow T^*M$. However, this is very awkward to work with.

We use the map $I: T(T^*M)\to T^{\bullet}(TM)$, see \cite[9.3.2]{Mackenzie:GT},
a \dvb isomorphism that induces the identity map on both side bundles and on the
core \vb.

Take the function
\begin{equation*}
\ell_{T(Y)}\circ I: T(T^*M)\rightarrow \R.
\end{equation*}
It follows directly that this is also linear with respect to $T^*M$.

Therefore, it will define a linear section $\mathfrak{Y}$ of the dual of the \vb
$T(T^*M)\rightarrow T^*M$, that is, of the \vb $T^*(T^*M)\rightarrow T^*M$.

Consider $\mathfrak{Y}(\ph)\in T^*(T^*M)$ for $\ph\in T^*M$. Pair this with a
$\xi\in T(T^*M)$ with outline $(\xi;\ph,x;m)$, where $x\in TM$. Using
(\ref{definition_of_eta_square_cap}), 
\begin{equation*}
\langle \mathfrak{Y}(\ph),\xi\rangle_{T^*M} = 
(\ell_{T(Y)}\circ I)(\xi) = 
\langle I(\xi),T(Y)(x)\rangle_{TM}.
\end{equation*}
By \cite[9.3.2]{Mackenzie:GT} and \cite[3.4.6]{Mackenzie:GT}, it follows that
\begin{equation}\label{eq:def-of-T(Y)}
\langle I(\xi),T(Y)(x)\rangle_{TM} 
=\langle d\ell_Y(\ph), \xi\rangle.
\end{equation}
This is true for any such $\xi\in T(T^*M)$ so it follows that
\begin{equation*}
\mathfrak{Y}(\ph) = (d\ell_Y)(\ph),
\end{equation*} 
and the linear section in question, $(T(Y)^{\sqcap},Y)$,
can be identified with $(d\ell_Y,Y)$.

Now $\wtilde{X}$ defines the linear function
\begin{equation*}
\ell_{\wtilde{X}}: T^*(TM)\rightarrow \R.
\end{equation*}
The function $\ell_{\wtilde{X}}$ is linear with respect to both $TM$ and
$T^*M$ (see Proposition~\ref{prop:Crelle}). 
Similarly as before, we elaborate on $\ell_{\wtilde{X}}$ being linear
with respect to $T^*M$. This defines a section of the dual of the \vb
$T^*(TM)\to T^*M$, that is, of $T^*(TM)\duer T^*M\to T^*M$. Again, this is not
easy to work with, and in this case we use the reversal isomorphism $R:
T^*(T^*M)\rightarrow T^*(TM)$.
It follows that
\begin{equation*}
\ell_{\wtilde{X}}\circ R: T^*(T^*M)\to \R
\end{equation*}
is also linear with respect to $T^*M$,
and defines a section $\mathfrak{X}$ of the dual of the \vb
$T^*(T^*M)\to T^*M$, that is, of the \vb $T(T^*M)\to T^*M$.

Then for $\ph \in T^*M$, and any $\mathfrak{F}\in T^*(T^*M)$ with outline
$(\mathfrak{F};x,\ph;m)$, with $x\in TM$, using
(\ref{definition_of_eta_square_cap}),
\begin{equation*}
\langle \X(\ph),\mathfrak{F}\rangle_{T^*M} = (\ell_{\wtilde{X}}\circ
R)(\mathfrak{F}) = \ell_{\wtilde{X}}(R(\mathfrak{F})) = \langle
R(\mathfrak{F}),\wtilde{X}(x)\rangle_{TM}.
\end{equation*}

At this point, we use that $R = J^*\circ I \circ (d\nu)^{\sharp}$,
\cite[p. 442]{Mackenzie:GT}, where $(d\nu)^{\sharp}$ is the map associated to
the canonical symplectic structure $d\nu$ on $T^*M$:
\begin{equation}
\label{eq:RJIdnu}  
\begin{tikzcd}[row sep=2cm, column sep = 2cm]
T^*(T^*M) \arrow[r,"R"]\arrow[d, swap,"(d\nu)^{\sharp}"] 
&T^*(TM) \arrow[d,<-,"J^*"]\\
T(T^*M)\arrow[r,swap,"I"]
&T^{\bullet}(TM)
\end{tikzcd}
\end{equation}
Note that (\ref{eq:RJIdnu}) is a commutative diagram of isomorphisms of \dvbs. 
Therefore, we can write
\begin{equation*}
\langle R(\mathfrak{F}),\wtilde{X}(x)\rangle_{TM} = 
\langle J^*(I ((d\nu)^{\sharp}(\mathfrak{F}))), 
\wtilde{X}(x)\rangle_{TM} = 
\langle I ((d\nu)^{\sharp}(\mathfrak{F})),  J(\wtilde{X}(x))\rangle_{TM}\\ = 
\langle I ((d\nu)^{\sharp}(\mathfrak{F})),  T(X)(x)\rangle_{TM}.
\end{equation*}
As before, using (\ref{eq:def-of-T(Y)}),
\begin{equation*}
\langle I ((d\nu)^{\sharp}(\mathfrak{F})),  T(X)(x)\rangle_{TM}=
 (d\nu)^{\sharp}(\mathfrak{F})(\ell_X) =
 \langle d\ell_X,(d\nu)^{\sharp}(\mathfrak{F}) \rangle = 
-\langle (d\nu)^{\sharp}(d\ell_X),\mathfrak{F}\rangle,
\end{equation*}
so we see that $\mathfrak{X} = -(d\nu)^{\sharp}(d\ell_X)$; that is, it is the
Hamiltonian vector field for the function $\ell_X$. Denote it by
$H_{\ell_X}$. Finally,
\begin{equation*}
\lpair T(Y)^{\sqcap}, \wtilde{X}^{\sqcap}\rpair_{T^*M} = \langle d\ell_Y,
H_{\ell_X}\rangle = \ell_{[X,Y]}.
\end{equation*}

\subsection{Example with $TA$}

In the case of Example \ref{ex:conn}, what are the corresponding
$T(\mu)^{\sqcap}$ and $(X^{\hlift})^{\sqcap}$ defined?

Just as in the case of $T^2M$, $(T(\mu)^{\sqcap},\mu)$ can be identified with
$(d\ell_{\mu},\mu)$.

We elaborate a bit more on $(X^{\hlift},X)$.
Again, we use $\ell_{X^{\hlift}}\circ R: T^*A^*\rightarrow \R$, and its
linearity with respect to $A^*$. This will define a section of the dual of
$T^*A^*\rightarrow A^*$, that is, of $TA^*\rightarrow A^*$. Denote this by
$\Phi$.

Then for $\kappa\in A^*$, $\Phi(\kappa)\in TA^*$. Pair it with a $\Psi\in
T^*A^*$ with outline $(\Psi;\kappa,a;m)$.

Then, for suitable $\mathcal{X}\in TA^*$, and using \cite[9.5.1]{Mackenzie:GT}, it
follows that
\begin{equation*}
\langle \Phi(\kappa),\Psi\rangle_{A^*} = 
(\ell_{X^{\hlift}}\circ R)(\Psi) = 
\langle R(\Psi), X^{\hlift}(a)\rangle_A = 
\llangle \mathcal{X}, X^{\hlift}(a)\rrangle - \langle \Psi,
\mathcal{X}\rangle_{A^*}.
\end{equation*}
The outlines for the elements involved are:
\begin{equation*}
\begin{tikzcd}[row sep=1cm, column sep = 1cm]
T^*A\ni R(\Psi) \arrow[r,mapsto," "]\arrow[d,mapsto, swap," ",xshift
= 4.5mm] &a\in A \arrow[d,mapsto," "]\\
A^*\ni \kappa\arrow[r,mapsto,swap," "]
&m, 
\end{tikzcd}\quad
\begin{tikzcd}[row sep=1cm, column sep = 1cm]
TA\ni X^{\hlift}(a)\arrow[r,mapsto," "]\arrow[d,mapsto, swap," ",xshift =
4.5mm] &X(m)\in TM \arrow[d,mapsto," "]\\
A\ni a\arrow[r,mapsto,swap," "]
&m, 
\end{tikzcd}
\end{equation*}
\begin{equation*}
\begin{tikzcd}[row sep=1cm, column sep = 1cm]
T^*A^*\ni \Psi \arrow[r,mapsto," "]\arrow[d,mapsto, swap," ",xshift =
6.5mm] &a\in A \arrow[d,mapsto," "]\\
A^*\ni \kappa\arrow[r,mapsto,swap," "]
&m,
\end{tikzcd}\quad
 \begin{tikzcd}[row sep=1cm, column sep = 1cm]
 TA^*\ni\mathcal{X} \arrow[r,mapsto," "]\arrow[d,mapsto, swap," ",xshift =
5mm] &X(m)\in TM \arrow[d,mapsto," "]\\
 A^*\ni \kappa\arrow[r,mapsto,swap," "]
 &m. 
 \end{tikzcd}
\end{equation*}
Now use \cite[3.4.6]{Mackenzie:GT} for the first term of the last equation.
Choose a $\ph\in \Gamma A^*$ with $\ph(m) = \kappa$, and a $\mu\in \Gamma A$
with $\mu(m) = a$. We can also make the following choice; linear vector fields
of a \vb $A$ are in bijective correspondence with linear vector fields on its
dual bundle $A^*$ (see \cite[3.4.5]{Mackenzie:GT}). Therefore, take
$\mathcal{X}$ to be $X^{\hlift_*}(\ph(m))$, where $X^{\hlift_*}$ is the corresponding
linear vector field to $X^{\hlift}$. Then we can write
\begin{multline}\label{eq:ex-conn-calculations}
\llangle X^{\hlift_*}(\ph(m)), X^{\hlift}(a)\rrangle - \langle \Psi,
X^{\hlift_*}(\ph(m))\rangle_{A^*} \\=
X^{\hlift_*}(\ph(m))(\ell_{\mu}) + X^{\hlift}(\mu(m))(\ell_{\ph})-X(m)(\langle \ph,
\mu\rangle) - \langle \Psi,X^{\hlift_*}(\ph(m))\rangle_{A^*}.
\end{multline}
At this point, we recall that for $\ph\in \Gamma A^*$, and for $\mu\in \Gamma
A$,
\begin{equation*}
X^{\hlift}(\ell_{\ph}) = \ell_{\nabla_X^{\dual}(\ph)}\in
C^{\infty}(E),\qquad
X^{\hlift_*}(\ell_{\mu}) = \ell_{\nabla_X(\mu)}\in C^{\infty}(E^*), 
\end{equation*}
and of course the relation between $\nabla$ and $\nabla^{\dual}$,
\begin{equation*}
\langle \nabla_X^{\dual}(\ph), \mu\rangle = 
X(\langle \ph, \mu\rangle) - \langle\ph, \nabla_X(\mu)\rangle,
\end{equation*}
and the latter equation can be rewritten as
\begin{equation*}
\ell_{\nabla^{\dual}_X(\ph)}\circ \mu = X(\langle \ph, \mu\rangle) -
\ell_{\nabla_X(\mu)}\circ \ph.
\end{equation*}
Returning to (\ref{eq:ex-conn-calculations}),
\begin{multline*}
\llangle X^{\hlift_*}(\ph(m)), X^{\hlift}(a)\rrangle - \langle \Psi,
X^{\hlift_*}(\ph(m))\rangle_{A^*} \\= 
\ell_{\nabla_X(\mu)}(\ph(m))  + \ell_{\nabla^{\dual}_X(\ph)}(\mu(m)) -
X(m)(\langle\ph, \mu\rangle) - \langle \Psi,X^{\hlift_*}(\ph(m))\rangle_{A^*} = 
-\langle \Psi,X^{\hlift_*}(\ph(m))\rangle_{A^*}.
\end{multline*}

Finally, the pairing between $(d\ell_{\mu},\mu)$ and $(X^{\hlift_*},X)$ is,
\begin{equation*}
\langle X^{\hlift_*},d\ell_\mu\rangle =X^{\hlift_*}(\ell_{\mu}) =
\ell_{\nabla_X(\mu)}.
\end{equation*}
  
\newcommand{\noopsort}[1]{} \newcommand{\singleletter}[1]{#1} \def\cprime{$'$}
  \def\cprime{$'$}

\end{document}